\documentclass[12pt,a4paper]{article}
\usepackage{amsmath}
\usepackage{amsfonts}
\usepackage{amssymb}
\usepackage{graphicx}%
\renewcommand{\baselinestretch}{1}
\numberwithin{equation}{section}
\newtheorem{theorem}{Theorem}[section]

\newtheorem{lemma}{Lemma}[section]

\newtheorem{proposition}{Proposition}[section]
\newtheorem{remark}{Remark}[section]

\newenvironment{proof}[1][Proof]{\noindent\textbf{#1.} }{\ \rule{0.5em}{0.5em}}

\begin{document}

\title{On similarity solutions for boundary layer flows with prescribed heat flux}
\author{BERNARD BRIGHI$\dag$  and JEAN-DAVID HOERNEL$\ddagger$ }
\date{}
\maketitle

\begin{center}
Universit\'e de Haute-Alsace, Laboratoire de Math\'ematiques et Applications
\vskip 0,1cm
4 rue des fr\`eres Lumi\`ere, 68093 MULHOUSE (France)
\vskip 1,5cm
\end{center}

\begin{abstract}
This paper is concerned with existence, uniqueness and behavior of 
the solutions of the autonomous third order nonlinear differential equation 
$f^{\prime\prime\prime}+\left(  m+2\right)  ff^{\prime\prime}-\left(
2m+1\right)  f^{\prime2}=0$ on $\mathbb{R}^+$ with the boundary conditions $f(0)=-\gamma$, 
$f^{\prime}\left(  \infty\right)  =0$ and $f^{\prime\prime}(0)=-1$. This 
problem arises when looking for similarity solutions for boundary layer flows with 
prescribed heat flux. To study solutions we use some direct approach as well as blowing-up 
coordinates to obtain a plane dynamical system.
\end{abstract}

\renewcommand{\baselinestretch}{1}
\footnotetext{AMS 2000 Subject Classification: 34B15, 34C11, 76D10.}
\footnotetext{Key words and phrases: Third order differential equations, boundary value problems, 
blowing-up coordinates, plane dynamical systems. }
\footnotetext{$\dag$ b.brighi@uha.fr $\ddagger$ j-d.hoernel@wanadoo.fr}

\section{Introduction}
We consider the following third order non-linear autonomous differential
equation found in \cite{Pop}%
\begin{equation}
f^{\prime\prime\prime}+\left(  m+2\right)  ff^{\prime\prime}-\left(
2m+1\right)  f^{\prime2}=0 \label{equation}%
\end{equation}
with the boundary conditions
\begin{equation}
f(0)=-\gamma, \label{cond01}%
\end{equation}%
\begin{equation}
f^{\prime}\left(  \infty\right)  =0, \label{cond02}%
\end{equation}%
\begin{equation}
f^{\prime\prime}(0)=-1 \label{cond03}%
\end{equation}
where $f'(\infty):=\underset{t\rightarrow\infty}{\lim} f'(t)$.

This equation gives the similarity solutions for free convection boundary-layer flows along a 
vertical permeable surface with prescribed surface heating and mass transfer rate. The solutions 
depend on two parameters: $m$, the power-law exponent and $\gamma$, the mass 
transfer parameter. The case $\gamma=0$ corresponds to an impermeable wall, $\gamma<0$ to a fluid suction and $\gamma>0$ to a fluid injection.
In the following we are investigating for existence and uniqueness of the solutions of 
(\ref{equation})-(\ref{cond03}) according to the values of $m$ and $\gamma$. We also 
gives some results about the boundedness and behavior of the solutions.

The problem involving similarity solutions with prescribed surface temperature leads to a similar 
equation with $f'(0)=1$ instead of $f''(0)=-1$ and is investigated in \cite{brighi02}, \cite{brighi01} and 
\cite{BrighiSari}. This alternative set of boundary conditions leads to significant differences in the obtained results and modelizes some very different physical problem (see \cite{Pop1} and \cite{Pop}
for more details about the physical 
interpretation of the two sets of boundary conditions). On the other hand the blowing-up coordinates introduced to transform the differential equation 
$(\ref{equation})$  are the same as in \cite{BrighiSari}, and the dynamical 
system obtained is very close to the one of \cite{BrighiSari}. For this reason we will refer to this paper for all that concerns the dynamical system.

The asymptotic behavior of the unbounded solutions for both prescribed surface temperature 
and prescribed heat flux is studied in \cite{BBequiv}.

\section{Preliminary results}

First, if $f$ verifies (\ref{equation}) let us notice that
\begin{equation}
\left(  f^{\prime\prime}e^{(m+2)F}\right)  ^{\prime}=\left(  2m+1\right)
f^{\prime2}e^{(m+2)F} \label{egalite01}%
\end{equation}
with $F$ any anti-derivative of $f.$ As $f^{\prime}$ and $f^{\prime\prime}$
cannot vanish at the same point without being identically equal to zero, we deduce the
\begin{lemma}
\label{concavite}Let $f$ be a non constant solution of $(\ref{equation})$ on some interval $I$. For all 
$t_{0}\in I$ we have
\begin{itemize}
\item If $m\leq-\frac{1}{2}$, $f^{\prime\prime}(t_{0})\leq0\Rightarrow f^{\prime\prime
}(t)<0$ for $t>t_{0}.$
\item If $m>-\frac{1}{2}$, $f^{\prime\prime}(t_{0})\geq0\Rightarrow f^{\prime\prime
}(t)>0$ for $t>t_{0}.$
\end{itemize}
\end{lemma}
\begin{proof}
It follows immediately from (\ref{egalite01}).
\end{proof}
\bigskip

Let us also remark that if $f$ is a solution of (\ref{ivp}) on $[0,T)$, then for $m\leq-\frac{1}{2}$ $f$
would be concave and for $m>-\frac{1}{2}$ it would be either concave or concave-convex.

\begin{proposition}
\label{prop<-1/2}For $m\leq-\frac{1}{2}$ there is only solutions to 
$(\ref{equation})$-$(\ref{cond03})$
if $f^{\prime}(0)>0.$ Moreover, if $f$ is a solution of $(\ref{equation})$-$(\ref{cond03})$ then
\begin{itemize}
\item $f$ is strictly concave, increasing and $f\left(  t\right)  \geq-\gamma$ for 
all $t$ in $[0,\infty)$.
\item If $m\in (-2,-\frac{1}{2}]$ and $\gamma>0$ then $f$ becomes positive at 
infinity. Moreover there exists
$t_{0}\geq\frac{\gamma}{f^{\prime}(0)}$ such that for all $t>t_{0},$ $f(t)>0.$
\end{itemize}
\end{proposition}
\begin{proof}
Since $f^{\prime\prime}(0)=-1$ and in view of lemma \ref{concavite}, $f^{\prime\prime}(t)$
would be negative for all $t$ which shows us that $f^{\prime}$
would be decreasing and $f$ concave. As we want to have $f^{\prime}\left(
\infty\right)  =0$, we must have $f^{\prime}(t)>0$ for all $t.$

For $m\in (-2,-\frac{1}{2}], $ using the fact that $f^{\prime\prime}(t)$
is negative for all $t,$ we see from (\ref{equation}) that if $f(t)\leq0$ for all $t$ we also 
have $f^{\prime \prime\prime}(t)<0$ for all $t.$ This implies that $f^{\prime}$ is concave 
and as $f^{\prime}$ is positive we cannot have $f^{\prime}(\infty)=0.$

Finally as $f$ is concave its graph is under its tangent in particular under
that at $0$ which equation is $y=f^{\prime}(0)t-\gamma.$ Thus $f$ becomes
positive after the point of intersection of its tangent at $0$ and the $t$-axis,
it means after $t_{0}=\frac{\gamma}{f^{\prime}(0)}.$
\end{proof}

\begin{proposition}
\label{prop>-1/2}Let $f$ be a solution of $(\ref{equation})$-$(\ref{cond03})$. For $m>-\frac{1}{2}$ 
we have
\begin{itemize}
\item Either $f$ is strictly concave and increasing and we must have $f^{\prime}\left(
0\right) >0$.
\item Or $f$ is concave-convex and

- if $f^{\prime}(0)\leq0$ the solution only exists for $\gamma<0$ and is
positive and decreasing.

- if $f^{\prime}(0)\geq0$ the solution is increasing-decreasing and positive
for $t\geq t_{0}$ with $t_{0}$ such that $f^{\prime\prime}(t_{0})=0.$
\end{itemize}
\end{proposition}
%\newpage
\begin{proof}
\begin{itemize}
\item As $f^{\prime\prime}(0)=-1$ if
$f^{\prime\prime}$ does not vanish it would remain negative and $f$
would be strictly concave. As above considering that $f^{\prime}$ would be decreasing,
 to have $f^{\prime}\left(  \infty\right)  =0$ we must
have $f^{\prime}>0$.
\item Suppose there exists $t_{0}$ such that $f^{\prime\prime}(t_{0})=0$ and
$f^{\prime\prime}<0$ on $\left[  0,t_{0}\right)  .$ Using lemma
\ref{concavite} we then have $f^{\prime\prime}>0$ on $\left(
t_{0},\infty\right)  $ which shows that $f$ would be concave-convex. We also
have that $f^{\prime}$ would be decreasing on $\left[  0,t_{0}\right)  $
and increasing on $\left[  t_{0},\infty\right)  $ which implies that
$f^{\prime}$ admits a negative minimum at $t_{0}$ because if not we cannot have
$f^{\prime}\left(  \infty\right)  =0.$ Thus we have the two following cases:
if $f^{\prime}(0)\leq0$ then $f^{\prime}(t)<0$ on $\left[  0,\infty\right)  $
and $f$ would be decreasing and if $f^{\prime}(0)\geq0$ then there
exists $t_{1}<t_{0}$ such that $f^{\prime}(t_{1})=0$ and $f$ would be
increasing on $\left[  0,t_{1}\right)  $ and decreasing on $\left[
t_{1},\infty\right)  $ which implies that $f$ admits a maximum at
$t_{1}.$ If now $f(t_{2})=0$ for some $t_{2} \geq t_{0}$, then 
$f(t)\leq0$  for all $t\geq t_{2}$ and since $f''(t)>0$ for $t>t_{2}$ we deduce
from (\ref{equation}) that $f'''\geq 0$ and that $f^{\prime}$ is convex on $[t_{2},\infty)$. 
But $f^{\prime}(t_{2})<0$
and we cannot have $f^{\prime}(\infty)=0,$ so $f(t)>0$ for all $t\geq
t_{0}.$ As a consequence we cannot have a concave-convex solution with
$f^{\prime}(0)\leq0$ and $\gamma>0.$
\end{itemize}
\vspace{-0.3cm}
\end{proof}

\begin{proposition}
For $m\geq -\frac{1}{2}$ the solutions of $(\ref{equation})$-$(\ref{cond03})$ are bounded.
\end{proposition}
\begin{proof}
For the concave-convex solutions the result is immediate. Suppose $f$ is
concave and unbounded, i.e. $f(t)\rightarrow \infty$ as $t \to \infty$. 
Then we have
\[%
\begin{array}
[c]{rl}
& f^{\prime\prime\prime}+\left(  m+2\right)  ff^{\prime\prime}=\left(
2m+1\right)  f^{\prime2}\geq0\\
\Rightarrow & f^{\prime\prime\prime}\geq-\left(  m+2\right)  ff^{\prime\prime}%
\end{array}
\]
and using the fact that $f^{\prime\prime}\leq0,$ if we choose $t_{1}$ such that
$f(t_{1})\geq \frac{1}{m+2}$ we have%
\begin{equation}
\forall t \in [t_{1},\infty), \quad f^{\prime\prime\prime}(t)\geq-f^{\prime\prime}(t). 
\label{>-1/2 bounded}%
\end{equation}
As $f^{\prime\prime\prime}\geq0$ on $\left[  t_{1},\infty\right)  ,$
$f^{\prime\prime}$ is increasing on $\left[  t_{1},\infty\right)  $ and using
the fact that $f^{\prime}\left(  \infty\right)  =0$ we deduce that
$f^{\prime\prime}\left(  \infty\right)  =0.$ Integrating $\left(
\ref{>-1/2 bounded}\right)  $ between the limits $r\geq t_{1}$ and $\infty$ leads to%
\[
\forall r\geq t_{1},\quad -f^{\prime\prime}(r)\geq f^{\prime}(r).
\]
Integrating once again we obtain
\[
\forall t\geq t_{1},\quad -f^{\prime}(t)+f^{\prime}(t_{1})\geq
f(t)-f(t_{1})
\]
which means that $f^{\prime}(\infty)=-\infty$ whereas one should have
$f^{\prime}\left(  \infty\right)  =0$, a contradiction.
\end{proof}

\begin{proposition} \label{f_seconde}
For all $m\in \mathbb{R}$ if $f$ is a solution of  $(\ref{equation})$-$(\ref{cond03})$ we have%
\[
\underset{t\rightarrow\infty}{\lim}f^{\prime\prime}(t)=0.
\]
\end{proposition}
\begin{proof}
See \cite{brighi01}.
\end{proof}

\begin{proposition}\label{blowT}
If a solution $f$ of $(\ref{equation})$ is only defined on a finite interval $[0,T)$, then $|f(t)|$, 
$|f'(t)|$ and $|f''(t)|$ tends toward infinity as $t \rightarrow T$.
\end{proposition}
\begin{proof}
See \cite{brighi01}.
\end{proof}

\subsection{Some equalities}

Integrating $\left(  \ref{equation}\right)  $ on $[\rho,r]$ leads
to%
\begin{equation}
f''(r)-f''(\rho)+(m+2)f(r)f'(r)-(m+2)f(\rho)f'(\rho)=3(m+1)\int_{\rho}^{r}f'(\xi)^2d\xi.
\label{int_01}%
\end{equation}
Multiplying $\left(  \ref{equation}\right)  $ by $t$ and integrating on
$[\rho,r]$ leads to%
\begin{align}
&  rf''(r)-\rho f''(\rho)-f'(r)+f'(\rho)+(m+2)(rf(r)f'(r)-\rho f(\rho)f'(\rho))
 \label{int_02}\\
& -\frac{(m+2)}{2}\left(  f(r)^2-f(\rho)^2\right)
 =3(m+1)  \int_\rho^r\xi f'(\xi)^2d\xi.\nonumber
\end{align}
Multiplying $\left(  \ref{equation}\right)  $ by $f$ and integrating on
$[\rho,r]$ leads to%
\begin{align}
&  f(r)f''(r)-f(\rho)f''(\rho)-\frac{1}{2}(f'(r)^2-f'(\rho)^2)
+(m+2)(f^2(r)f'(r)-f^2(\rho)f'(\rho))\label{int_03}\\
&  =(4m+5)  \int_{\rho}^{r}f(\xi)f'(\xi)^2d\xi.\nonumber
\end{align}

\subsection{The plane dynamical system}

Consider a right maximal interval $I=\left[  \tau,\tau+T\right)$ on which
$f$ does not vanish. For all $t$ in $I$, set
\begin{equation}
s=\int_{\tau}^{t}f(\xi)d\xi,\quad u(s)=\frac{f^{\prime}(t)}{f(t)^2}, \quad
v(s)=\frac{f^{\prime\prime}(t)}{f(t)^3},
\label{new_function}%
\end{equation}
to obtain the system
\begin{equation}
\left\{
\begin{array}
[c]{l}%
\dot{u}=P(u,v):=v-2u^{2},\\
\dot{v}=Q_{m}(u,v):=-(m+2)v+(2m+1)u^{2}-3uv,
\end{array}
\right.  \label{system}%
\end{equation}
in which the dot denotes the differentiation with respect to $s.$ Let us notice that if $f$ is negative on $I$ then $s$ decreases as $t$ grows.

The singular points of $\left(  \ref{system}\right)  $ are $O=(0,0)$ and
$A=\left(  -\frac{1}{2},\frac{1}{2}\right)  .$ The isoclinic curves $P(u,v)=0$ and
$Q_{m}(u,v)=0$ are given by $v=2u^{2}$ and $v=\Psi_{m}(u)$ where%
\[
\Psi_{m}(u)=\dfrac{(2m+1)u^{2}}{3u+(m+2)}.
\]

The point $A$ is

\begin{itemize}
\item An unstable node for $m\leq\frac{3-2\sqrt{6}}{2}$ ($\lambda_{1}\geq0$
and $\lambda_{2}\geq0$)$.$

\item An unstable focus if $\frac{3-2\sqrt{6}}{2}<m<\frac{3}{2}$
($Re(\lambda_{1})\geq0$ and $Re(\lambda_{2})\geq0$)$.$

\item A center if $m=\frac{3}{2}.$

\item A stable focus if $\frac{3}{2}<m<\frac{3+2\sqrt{6}}{2}$ ($Re(\lambda
_{1})\leq0$ and $Re(\lambda_{2})\leq0$)$.$

\item A stable node if $m\geq\frac{3+2\sqrt{6}}{2}$ ($\lambda_{1}\leq0$ and
$\lambda_{2}\leq0$)$.$\bigskip
\end{itemize}

For $m\neq-2,$ the singular point
$O$ is a saddle-node of multiplicity 2. It admits a center manifold $\mathcal W_{0}$ that is
tangent to the subspace $L_{0}=Sp\left\{  (1,0)\right\}$ and a stable (resp. unstable) 
manifold $\mathcal W$ if $m>-2$ (resp $m<-2)$ that is tangent to the subspace 
$L=Sp\left\{  (1,-(m+2))\right\}.$

\bigskip

We will now precise the phase portrait of the vector field $\left(
\ref{system}\right)  $ near the saddle-node point $O$ using the same arguments as in
\cite{BrighiSari} (see Fig 2.2.1).

\begin{itemize}
\item The parabolic sector is delimited by the separatrices $S_{0}$ and
$S_{1}$ which are tangent to $L$ at $O$.

\item The first hyperbolic sector is delimited by $S_{0}$ and the separatrix
$S_{2}$ which is tangent to $L_{0\text{ }}$at $O.$ The second hyperbolic
sector is delimited by $S_{1}$ and $S_{2}$.

\item The manifold $\mathcal W$ is the union of the separatrices $S_{0},$ $S_{1}$ and
the singular point $O$%
\[
\mathcal W=\left\{  S_{0}\right\}  \cup\left\{  O\right\}  \cup\left\{  S_{1}\right\}.\]
Near $O$, the manifold $\mathcal W$ takes place below $L$ for $m<-2$
or $m>-1$ and above $L$ for $-2<m<-1$.

In the case $m=-1$ the manifold $\mathcal W$ is given by
${\mathcal W}=\left\{  \left(  u,-u\right)  \in\mathbb{R}^{2};u>-\frac{1}{2}\right\}.$

\item The manifold $\mathcal W_{0}$ is the union of the separatrix $S_{2},$ the point
$O$ and a phase curve $C_{3}$%
\[
\mathcal W_{0}=\left\{  S_{2}\right\}  \cup\left\{  O\right\}  \cup\left\{
C_{3}\right\}  .
\]
Near the point $O$, the center manifold $\mathcal W_{0}$ takes place above $L_{0}$ for
$m<-2$ or $m>-\frac{1}{2}$ and below $L_{0}$ for $-2<m<-\frac{1}{2}$.

For $m=-\frac{1}{2},$ the center manifold $\mathcal W_{0}$ coincides with the $u$-axis.
\end{itemize}

\begin{remark}
We will not consider the case $m=-2$ because we will see later that there is no solution.
\end{remark}

If we note $S_{i}^{+}$ for an $\omega$-separatrix and $S_{i}^{-}$ for an $\alpha$-separatrix 
the behavior of the vector field in the neighborhood of the saddle-node point $O$ 
is given by the following figures
\[
\begin{array}
[c]{cc}
\begin{array}
[c]{c}
\raisebox{-0cm}{\includegraphics[scale=.3]
{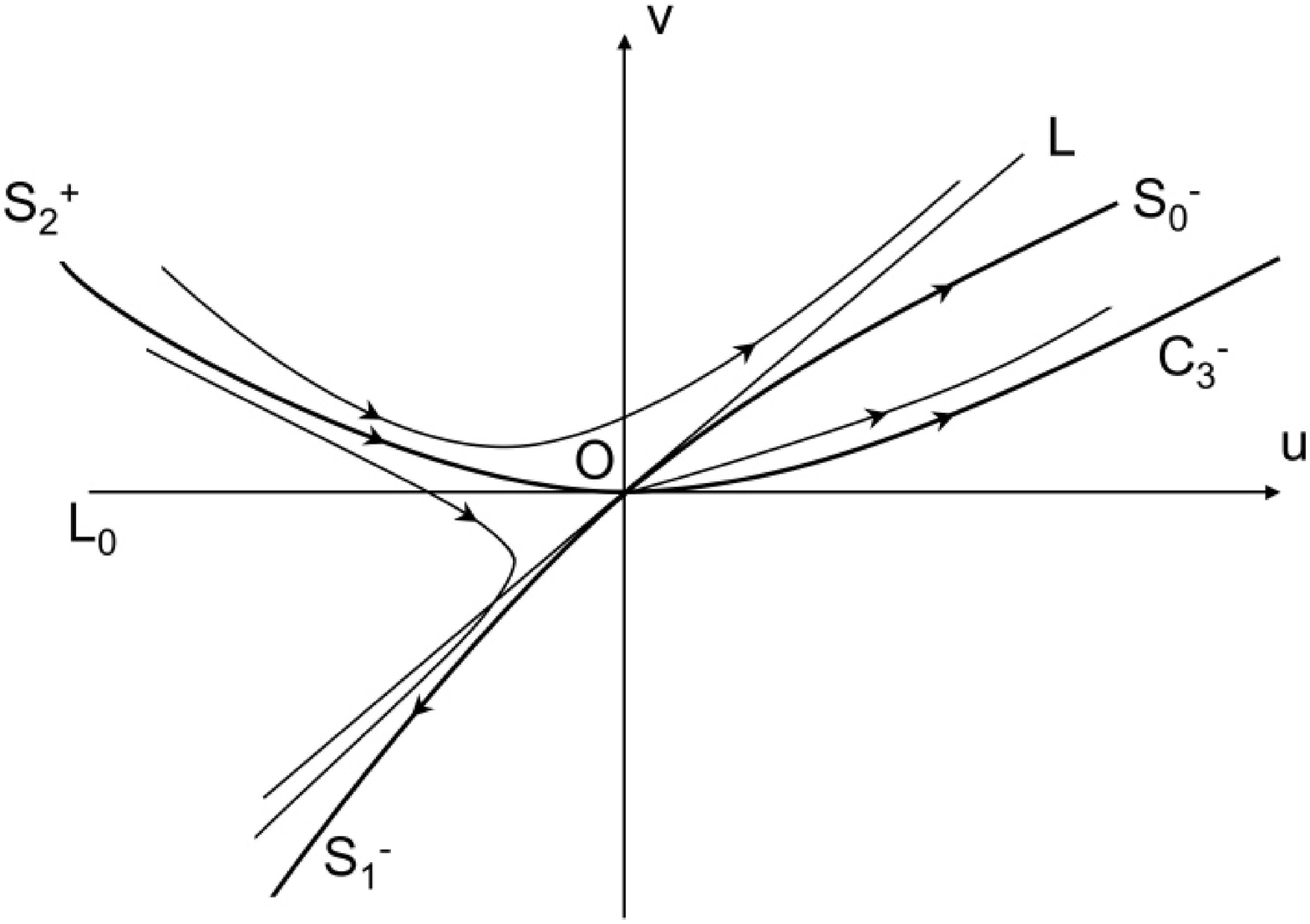}
}
\\
m<-2
\end{array}
&
\begin{array}
[c]{c}
\raisebox{-0cm}{\includegraphics[scale=.3]
{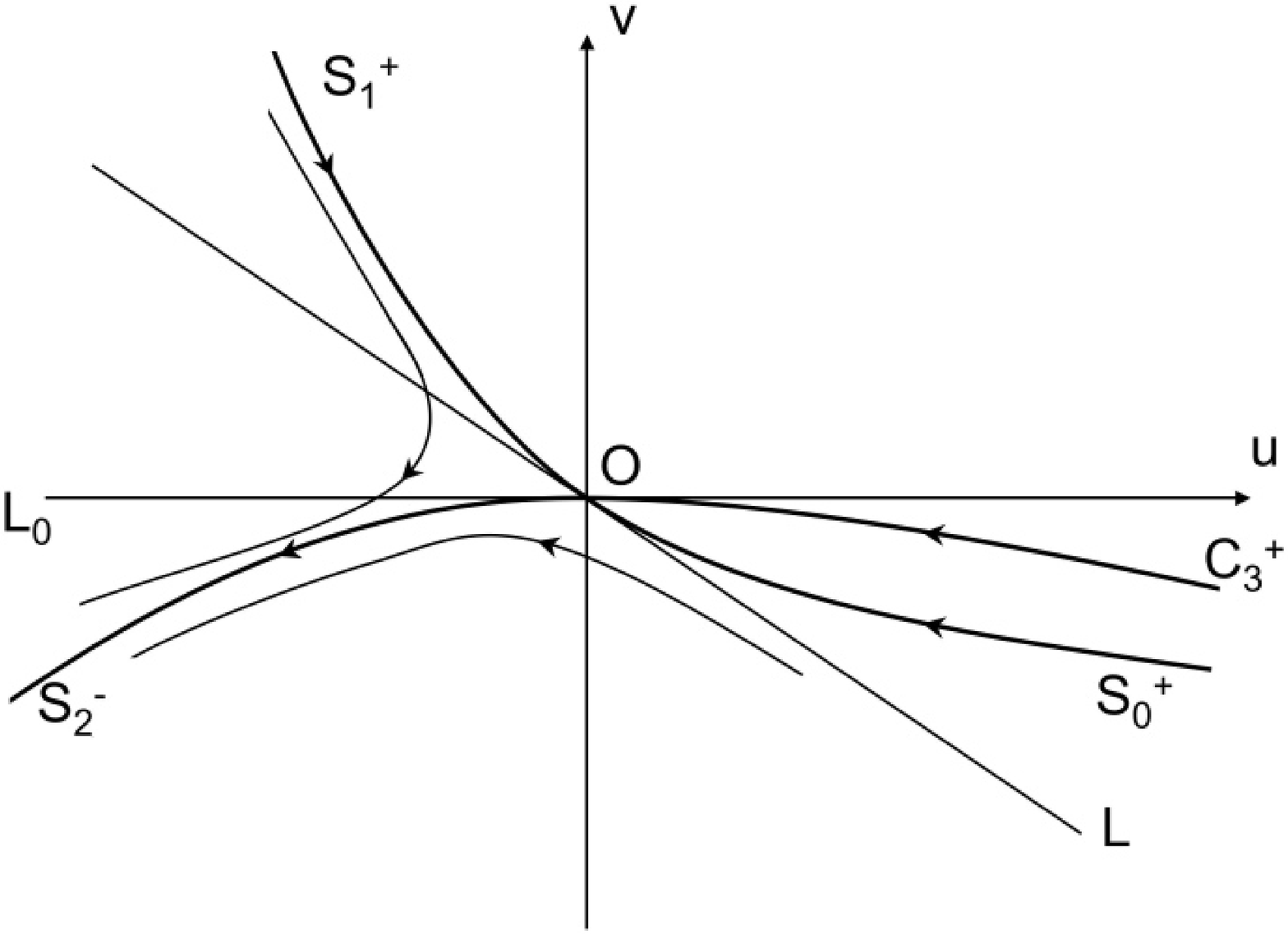}
}
\\
-2<m<-1
\end{array}
\\
\end{array}
\]

\[
\begin{array}
[c]{cc}
\begin{array}
[c]{c}
\raisebox{-0cm}{\includegraphics[scale=.3]
{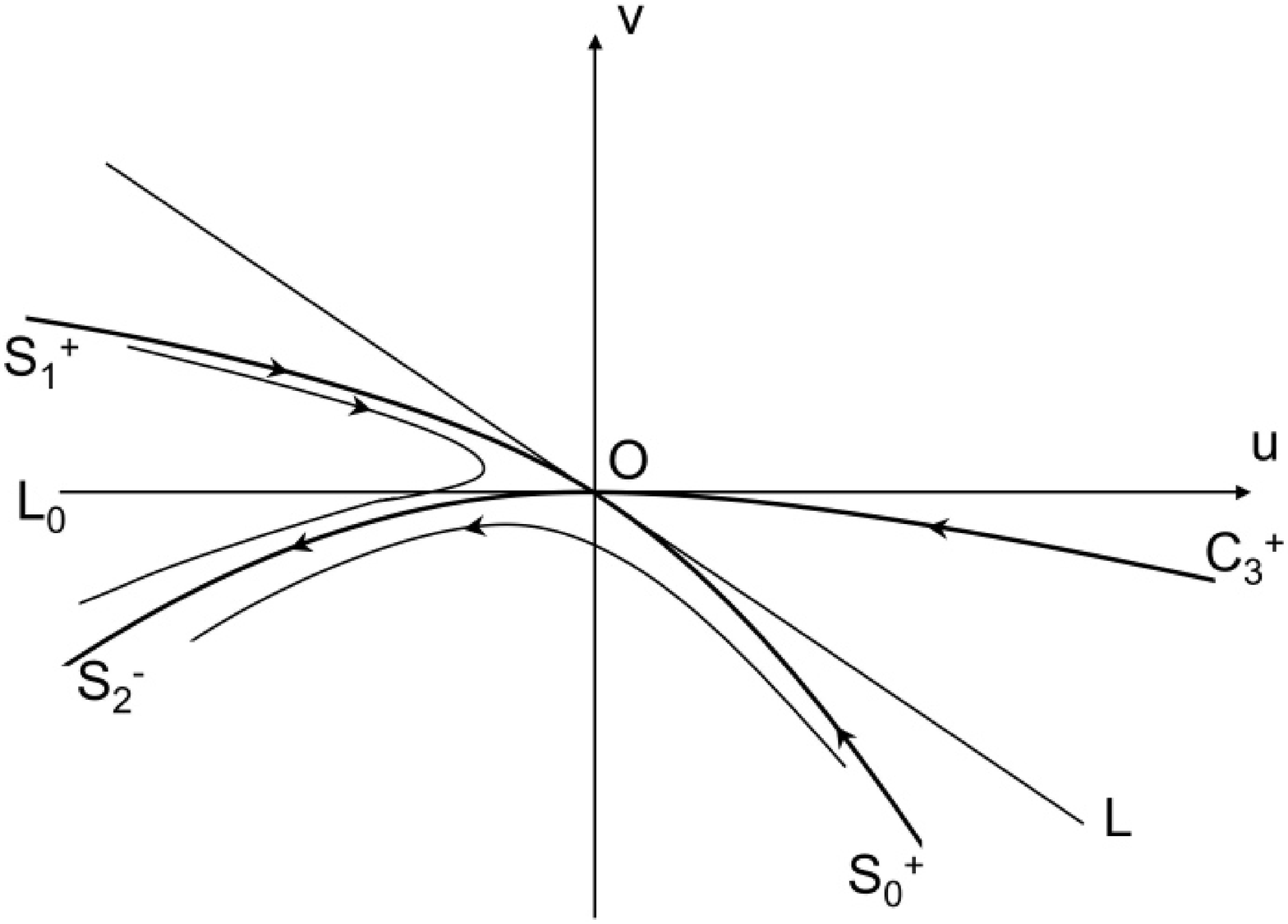}
}
\\
-1<m<-\frac{1}{2}
\end{array}
&
\begin{array}
[c]{c}
\raisebox{-0cm}{\includegraphics[scale=.3]
{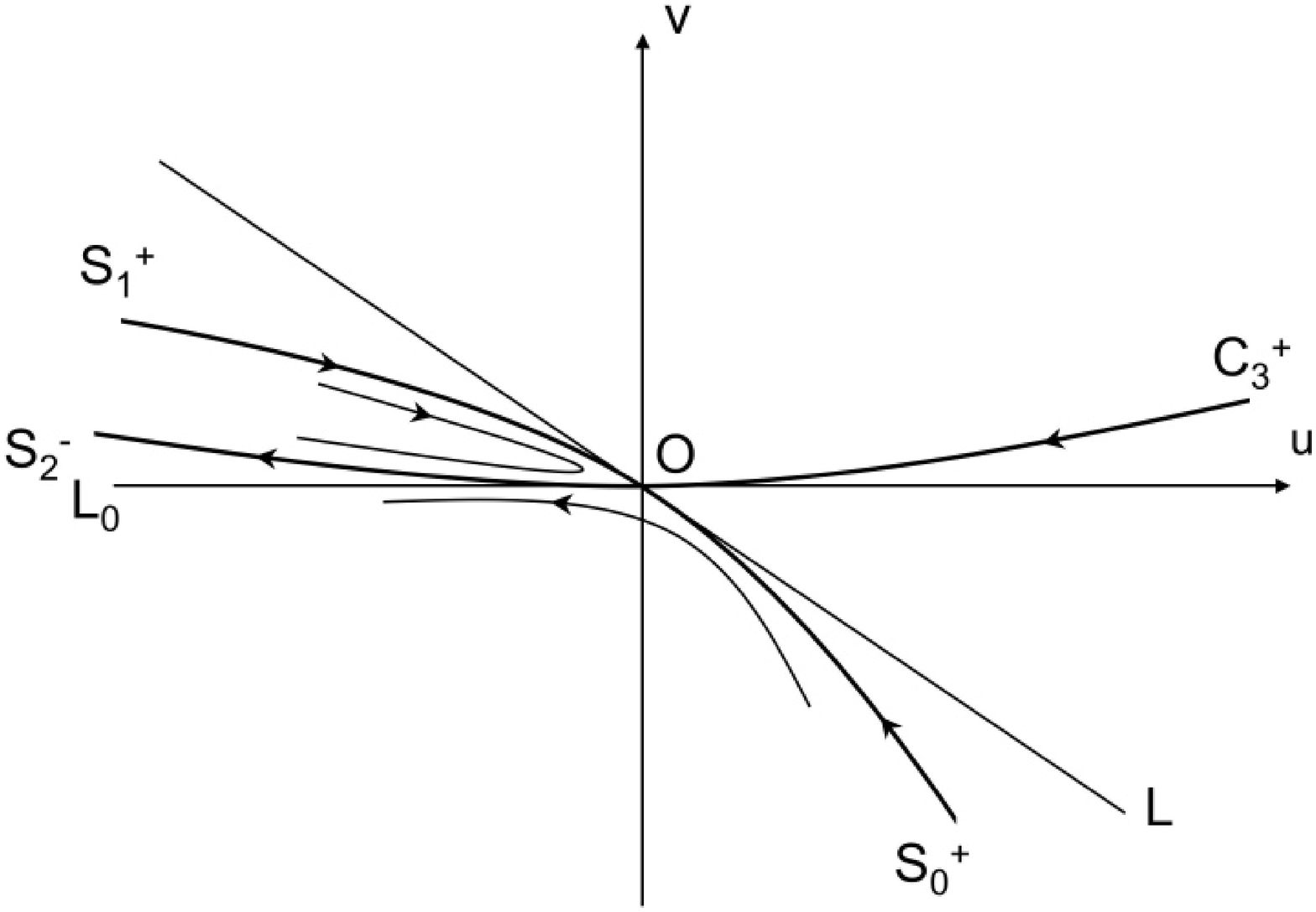}
}
\\
m>-\frac{1}{2}
\end{array}
\end{array}
\]
\centerline{Fig 2.2.1}

\newpage
We will also use the following notations: Consider a connected piece of a phase
curve $C$ of $\left(  \ref{system}\right)  $ lying in the region $P(u,v)<0$
(resp. $P(u,v)>0$), then $C$ can be characterized by $v=V_{m}(u)$ (resp.
$v=W_{m}(u)$) with $u$ belonging to some interval and $V_{m}$ (resp.
$W_{m}$) a solution of the differential equation
\begin{equation}
v^{\prime}=F_{m}(u,v):=\frac{Q_{m}(u,v)}{P(u,v)}=\frac{-(m+2)v+(2m+1)u^{2}%
-3uv}{v-2u^{2}}. \label{Fm}%
\end{equation}
\bigskip

\section{Main results}

To obtain results about the boundary value problem (\ref{equation})-(\ref{cond03}) we will often 
use the initial value problem $\mathcal{P}_{m,\gamma,\alpha}$
\begin{equation}%
\left \{
\begin{array}
[c]{rl}%
&f'''+(m+2)ff''-(2m+1)f'^2 =0,\\
&f(0) =-\gamma,\\
&f'(0) =\alpha,\\
&f''(0) =-1
\end{array}
\right .
\label{ivp}%
\end{equation}
with $\alpha \in \mathbb{R}$.

\subsection{The case $m\leq-2$}

\begin{lemma}\label{m<-2}
Let $m\leq-2$. If $\gamma\leq\sqrt[3]{\frac{2}{\left(  m+2\right)  ^{2}}}$ the
problem $(\ref{equation})$-$(\ref{cond03})$ has
no solution. In particular, for $m=-2$ there is no solution at all.
Moreover, if $\gamma>\sqrt[3]{\frac{2}{\left(  m+2\right)  ^{2}}}$ and if $f$ is a solution 
of $(\ref{equation})$-$(\ref{cond03})$, then $f$ is negative 
and the phase curve $(u(s),v(s))$ defined by $\left(  \ref{new_function}\right)  $ with
$\tau=0$ is a negative semi-trajectory which lies for $-s$ large enough in the
bounded domain
\[
\mathcal{D}_{+}=\left\{  \left(  u,v\right)  \in\mathbb{R}^{2}\, ;\quad  0<u<-\frac
{m+2}{2}\quad \text{and}\quad 0\leq v<-(m+2)u\right\}  .
\]
\end{lemma}
\begin{proof}
Suppose that $f$ is a solution of (\ref{equation})-(\ref{cond03}). By proposition
\ref{prop<-1/2}, we know that $f''(t)<0$ and $f'(t)>0$ for all $t$. If there exists $t_{1}$
such that $f(t)>0$ for $t>t_{1}$, then we deduce from (\ref{equation}) that $f'''(t)\leq 0$
for $t>t_{1}$. This implies that $f'$ is concave on $(t_{1},\infty)$, which does not allow to
have $f'(t_{1})>0$ and $f'(\infty)=0$. Therefore, if $f$ is a solution of 
(\ref{equation})-(\ref{cond03}) we necessarily have $\gamma>0$ and $f(t)<0$ for all $t$.
Next, we have 
\begin{equation}
\forall t\geq0,\quad \frac{f^{\prime}(t)}{f(t)^{2}}\geq0\quad \text{and}\quad %
\frac{f^{\prime\prime}(t)}{f(t)^{3}}>0. \label{l0}%
\end{equation}
As $f$ is bounded, we can write $(\ref{int_01})$ with $\rho=t$ and $r=\infty$ to get
\[
\forall t\geq0,\quad f^{\prime\prime}(t)+(m+2)f(t)f^{\prime}(t)=-3(m+1)\int
_{t}^{\infty}f^{\prime}(\xi)^{2}d\xi,
\]
and
\begin{equation}
\forall t\geq0,\quad f^{\prime\prime}(t)+(m+2)f(t)f^{\prime}(t)>0
\label{l1}%
\end{equation}
as $m<-2$. Let $\lambda$ be the limit of $f$ at infinity, integrating again leads to%
\begin{equation}
\forall t\geq0,\quad f^{\prime}(t)+\frac{(m+2)}{2}f(t)^{2}<\frac{(m+2)}%
{2}\lambda^{2}<0. \label{l2}%
\end{equation}
Writing $\left(  \ref{l1}\right)  $ and $(\ref{l2})$ for $t=0$ we obtain%
\[
\gamma^{2}>-\frac{2f'(0)}{m+2}\quad \text{and}\quad f'(0)>-\frac{1}{(m+2)\gamma},
\]
and finally $\gamma>\sqrt[3]{\frac{2}{\left(  m+2\right)  ^{2}}}$. Dividing
$\left(  \ref{l1}\right)  $ by $f(t)^{3}$ and $(\ref{l2})$ by $f(t)^{2}$ we get
\begin{equation}
\forall t\geq0,\quad \frac{f^{\prime}(t)}{f(t)^{2}}+\frac{(m+2)}{2}<0
\quad \text{and}\quad 
\frac{f^{\prime\prime}(t)}{f(t)^{3}}+(m+2)\frac{f^{\prime}(t)}{f(t)^{2}}<0. \label{l3}%
\end{equation}
Using the first of the two precedent inequalities we found that%
\[
\forall t\geq0,\quad f(t)\leq\frac{1}{\frac{m+2}{2}t-\frac{1}{\gamma}}%
\]
which implies that $\int_{0}^{\infty}f(\xi)d\xi=-\infty.$ Hence the trajectory
$s\mapsto(u(s),v(s))$ is defined on the whole interval $\left( -\infty
,0\right]  $ and with $(\ref{l0})$ and $(\ref{l3})$ the proof is complete.
\end{proof}
\bigskip

In the following we will sometimes need the system (\ref{system}) to obtain results 
about the problem (\ref{equation})-(\ref{cond03}) when direct approach fails. To this end we will give the behavior of the separatrices without proof because it is the same as in \cite{BrighiSari}.

\begin{theorem}
Let $m<-2.$ There exists $\gamma_{\ast}$ such that the problem 
$(\ref{equation})$-$(\ref{cond03})$ has infinitely many
solutions if $\gamma>\gamma_{\ast}$, one and only one solution if
$\gamma=\gamma_{\ast}$ and no solution if $\gamma<\gamma_{\ast}.$
\end{theorem}
\begin{proof}
Taking into account proposition \ref{prop<-1/2} and lemma \ref{m<-2}, consider
the solution of the initial value problem $(\ref{ivp})$ with $\alpha>0$ and $\gamma>0$.
Denote by $C_{\gamma,\alpha}$ the corresponding trajectory of the
plane system $\left(  \ref{system}\right)  $ defined by $\left(
\ref{new_function}\right)  $ with $\tau=0.$ We have $u(0)=\frac{\alpha}{\gamma^2}>0$
and $v(0)=\frac{1}{\gamma^3}>0$.

Before going further, according to Fig 2.2.1 we just have to precise the behavior of the separatrice $S_0^-$. As $s$ grows from $-\infty$, the $\alpha$-separatrix $S_{0}^{-}$ leaves to the
right the singular point $O$ tangentially to $L$ and intersects first the isocline
$Q_{m}(u,v)=0,$ then the isocline $P(u,v)=0$, the $u$-axis and the $v$-axis (see Fig 3.1.1).
Let $\left(  u_{\ast},v_{\ast}\right)  $ 
be the point where the separatrix $S_{0}^{-}$ intersects the isocline
$Q_{m}(u,v)=0$ and set $\gamma_{\ast}=\frac{1}%
{\sqrt[3]{v_{\ast}}}.$

If $\gamma<\gamma_{\ast},$ the straight line $v=\frac{1}{\gamma^{3}}$ does not
intersect the separatrix $S_{0}^{-}$ and for all $\alpha>0$, the Poincar\'e-Bendixson theorem
shows that $C_{\gamma,\alpha}$ does not remain in the bounded domain $\mathcal{D}_{+}$.
It follows from lemma \ref{m<-2} that $f$ cannot be a solution of 
(\ref{equation})-(\ref{cond03}) for any $\alpha>0$.

For $\gamma=\gamma_{\ast}$ the straight line $v=\frac{1}{\gamma_{\ast}^{3}}$
intersects the separatrix $S_{0}^{-}$ at the point $\left(  u_{\ast},v_{\ast
}\right)  $. As above, $f$ is not a solution for $\alpha\neq\gamma^{2}u_{\ast}$. For
$\alpha=\gamma^{2}u_{\ast},$ the phase curve $C_{\gamma,\alpha}$ is a negative
semi-trajectory which coincides with the part of the separatrix $S_{0}^{-}$ coming 
from $O$. Then $f$ cannot vanish, because on the contrary one of the coordinates $u$ or $v$
should go to infinity (recall $f'$ and $f''$ cannot vanish at the same point). Hence as long as $f$ exists we have $f'>0$ and $f''<0$, which implies that $f$ exists on the whole interval $[0,\infty)$.
Moreover $f^{\prime}(t)\rightarrow l\geq0$ as
$t\rightarrow\infty$ and supposing $l>0$ leads to a contradiction due to the negativity of
$f$. Therefore $f$ is a solution of $\left(\ref{equation}\right)  $-$\left(  \ref{cond03}\right)$.

%\begin{figure}
%\centering
%\includegraphics[natheight=5.76in,natwidth=7.68in,height=8.1166cm,width=10.7989cm]
%{m_inf_-2.jpg}
%\caption{$m<-2$}
%\end{figure} 
         
\[
%TCIMACRO{\FRAME{itbpFU}{10.7989cm}{8.1166cm}{0cm}{\Qcb{$m<-2$}}{}%
%{Figure}{\special{ language "Scientific Word";  type "GRAPHIC";
%maintain-aspect-ratio TRUE;  display "USEDEF";  valid_file "T";
%width 10.7989cm;  height 8.1166cm;  depth 0cm;  original-width 7.6803in;
%original-height 5.76in;  cropleft "0";  croptop "1";  cropright "1";
%cropbottom "0";  tempfilename 'm_inf_-2.jpg';tempfile-properties "XPR";}}}%
%BeginExpansion
\raisebox{-0cm}{\parbox[b]{10.7989cm}{\begin{center}
\includegraphics[scale=.4]
{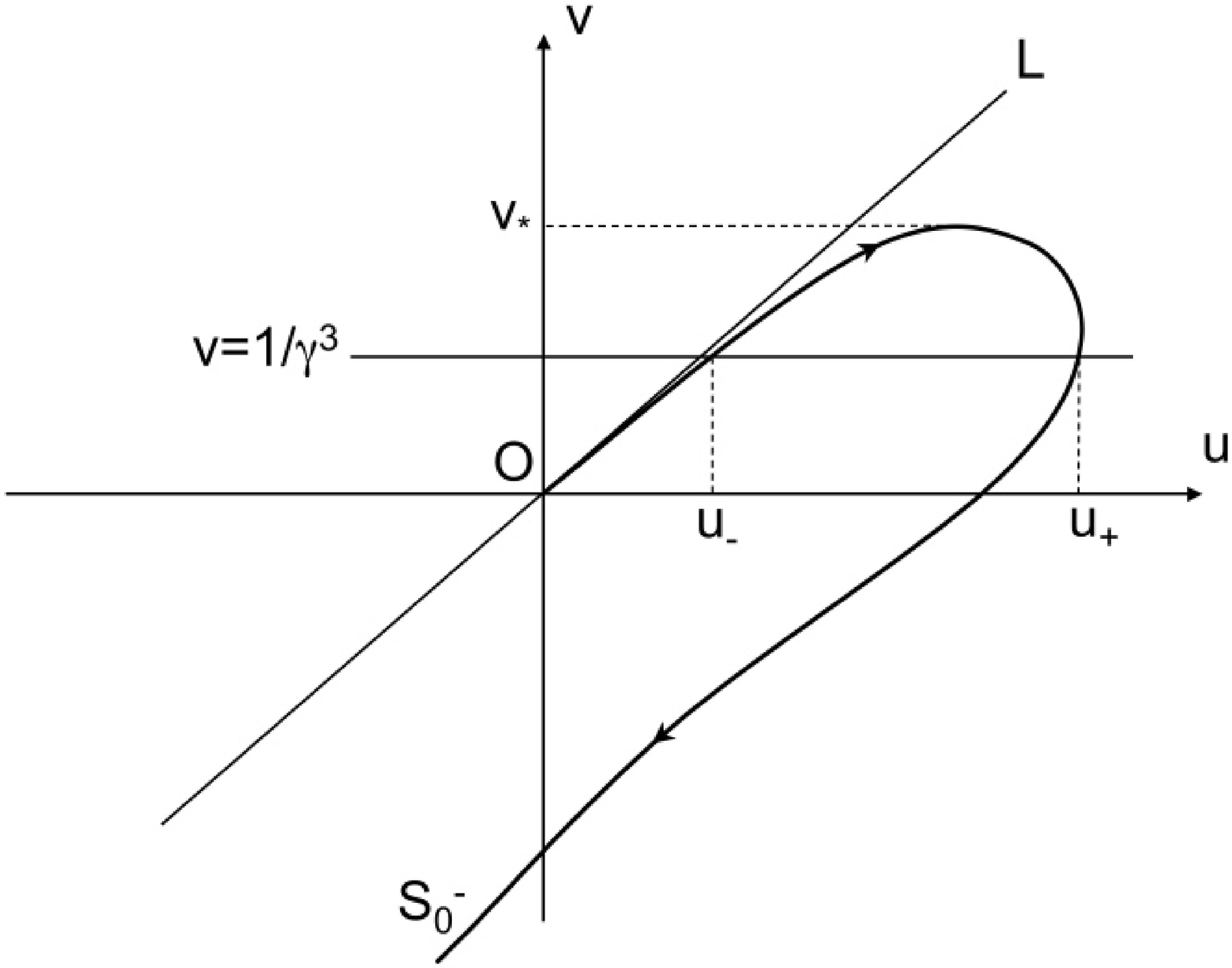}
\\
$m<-2$\\
{Fig 3.1.1}
\end{center}}}%
%EndExpansion
%\tag{($m<-2$)}%
\]

For $\gamma>\gamma_{\ast},$ the straight line intersects two times the
separatrix $S_{0}^{-}$ in $\left(  u_{-},\frac{1}{\gamma^{3}}\right)  $ and
$\left(  u_{+},\frac{1}{\gamma^{3}}\right)$ as shown in Fig 3.1.1. Using the same arguments as
above we see that if $\alpha\in\left[  u_{-}\gamma^{2},u_{+}\gamma^{2}\right]
$ then $f$ is a solution of $\left(  \ref{equation}\right)  $-$\left(
\ref{cond03}\right)  $ and if $\alpha\notin\left[  u_{-}\gamma^{2},u_{+}%
\gamma^{2}\right]  $ then $f$ is not. 
\end{proof}

\begin{remark}
Notice that for all $\gamma>0$ we have
$0<u_{+}\leq\frac{1}{\sqrt{2\gamma^{3}}}$ $($i.e. $0<\alpha\leq
\sqrt{\frac{\gamma}{2}})$ and that $u_{-}\rightarrow0$ as $\gamma\rightarrow\infty$.
\end{remark}

\begin{proposition}\label{limf}
Let $m<-2$ and $f$ be a solution of $(\ref{equation})$-$(\ref{cond03})$, then for 
$\gamma=\gamma_{*}$ we have that $f(t)\rightarrow \lambda <0$
as $t\rightarrow \infty$ and for every $\gamma> \gamma_{*}$ there are two solutions $f$
such that $f(t)\rightarrow \lambda <0$ as $t\rightarrow \infty$ and all the other solutions verify 
$f(t)\rightarrow 0$ as $t\rightarrow \infty$
\end{proposition}
\begin{proof}
The proof is the same as in \cite{BrighiSari}.
\end{proof}
\bigskip

\subsection{The case $-2<m<-1$}

\begin{proposition}\label{-2m-1}
For $-2<m<-1$ and $\gamma\geq0$ the problem $(\ref{equation})$-$(\ref{cond03})$ has no 
solution. Moreover, to have
solutions with $\gamma<0$ we must have $f^{\prime}(0)\geq-\frac{1}{(m+2)  \gamma}.$
\end{proposition}
\begin{proof}
From proposition \ref{prop<-1/2} if $f$ is a solution of $(\ref{equation})$-$(\ref{cond03})$
we know that $f'(0)>0$ and that for $t$ large enough
$f(t)>0$ and $f'(t)>0$. Thus for $-2<m<-1$ and $\gamma \geq 0$, we get from (\ref{int_01})
with $\rho=0$ and $r=t$
\begin{equation}
f^{\prime\prime}(t)=3\left(  m+1\right)  \int_{0}^{t}f^{\prime2}(\xi)d\xi
-\left(  m+2\right)  f(t)f^{\prime}(t)-\left(  m+2\right)
\gamma f^{\prime}(0)-1\leq -1.
\label{int_05}
\end{equation}
and a contradiction with proposition \ref{f_seconde}.

Let $\gamma<0$, then for all $t>0$ we have $f'(t)>0$ and $f(t)>0$. Using equality (\ref{int_05}) 
we obtain
$$-f''(t)-(m+2)f'(0)\gamma-1\geq 0$$
and as $t$ goes to infinity, using proposition \ref{f_seconde} this inequality gives
the second part of the result.
\end{proof}

\begin{theorem} \label{m_inf_-5_4} For $-2<m<-1$ there exists $\gamma_{*}<0$ 
such that the problem $(\ref{equation})$-$(\ref{cond03})$ has no solutions for $\gamma>\gamma_{*}$,
one and only one solution which is bounded for $\gamma=\gamma_{*}$ and two bounded solutions 
and infinitely many unbounded solutions for $\gamma<\gamma_{*}$.
\end{theorem}
\begin{proof}
From proposition  \ref{-2m-1} we know that if $\gamma\geq0$ there is no solution, 
so we must consider a solution $f$ of the initial value problem (\ref{ivp}) with $\gamma<0$ and 
$\alpha>0$. 
Let $C_{\gamma,\alpha}$ be the phase curve corresponding to $u,$ $v$ defined by (\ref{new_function}) with $\tau=0$. Looking at Fig 2.2.1, we 
see that the separatrix $S_{0}^+$ crosses first the $u$-axis, then the isocline $Q_{m}(u,v)=0$
before going to $O$. 
Let us call $(u_{*},v_{*})$ the point where the separatrix $S_{0}^+$ crosses the isocline $Q_{m}(u,v)=0$ and set $\gamma_{*}=\frac{1}{\sqrt[3]{v_{*}}}$ (see Fig 3.2.1).

\[%
%TCIMACRO{\FRAME{itbpFU}{10.7968cm}{8.1166cm}{0cm}{\Qcb{$-2<m<-1$}}%
%{}{Figure}{\special{ language "Scientific Word";  type "GRAPHIC";
%maintain-aspect-ratio TRUE;  display "USEDEF";  valid_file "T";
%width 10.7968cm;  height 8.1166cm;  depth 0cm;  original-width 9.5998in;
%original-height 7.1996in;  cropleft "0";  croptop "1";  cropright "1";
%cropbottom "0";  tempfilename 'm_inf_-1.jpg';tempfile-properties "XPR";}}}%
%BeginExpansion
\raisebox{-0cm}{\parbox[b]{10.7968cm}{\begin{center}
\includegraphics[scale=.4]
{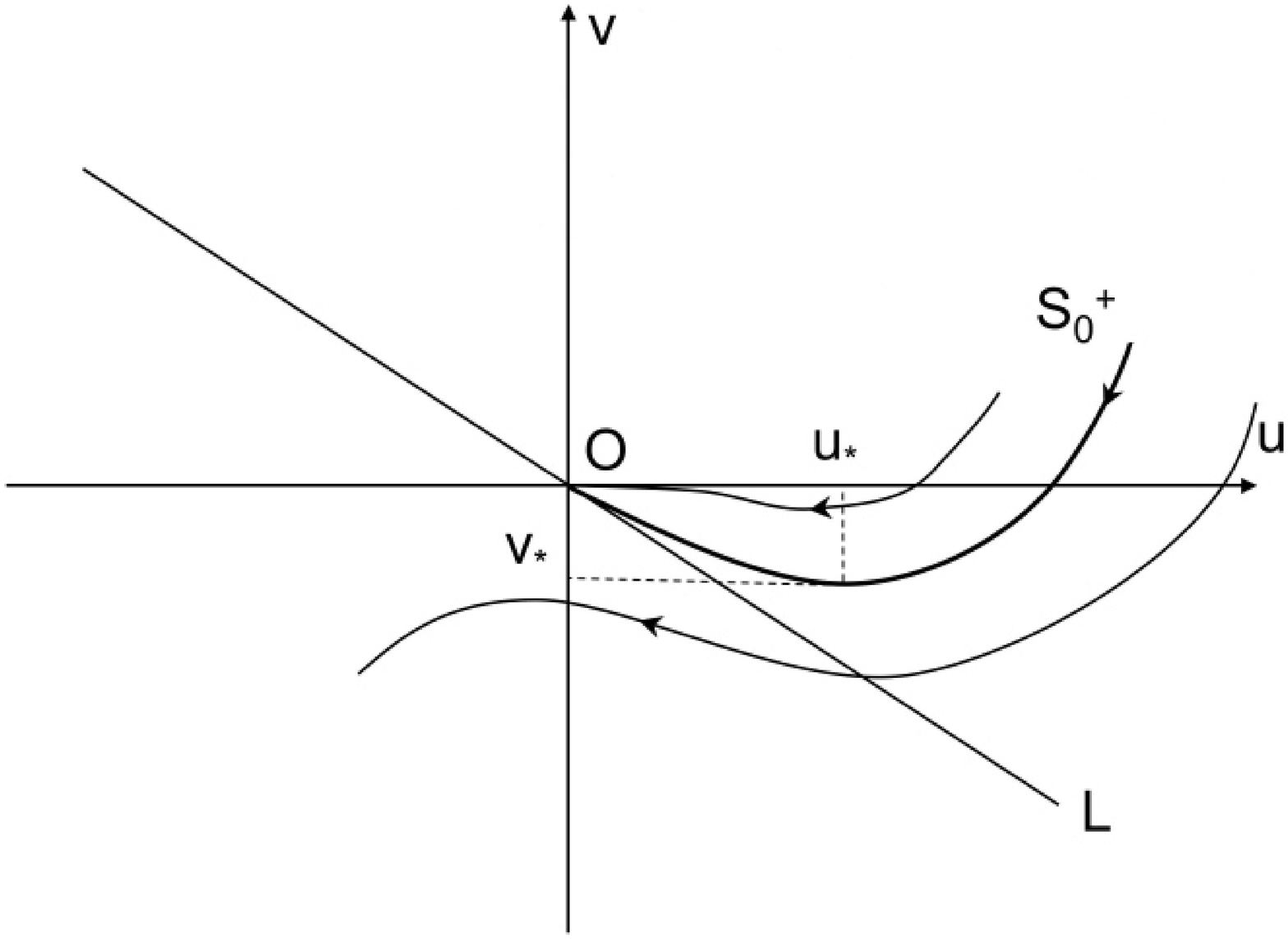}
\\
$-2<m<-1$\\
{Fig 3.2.1}
\end{center}}}%
%EndExpansion
\]
We see that the horizontal line $v=\frac{1}{\gamma^3}$ does not intersect the separatrix $S_{0}^+$
if $\gamma_{*}<\gamma<0$, is tangent to it if $\gamma=\gamma_{*}$, and intersects it through two 
points $(u_{-},\frac{1}{\gamma^3})$ and $(u_{+},\frac{1}{\gamma^3})$ if $\gamma<\gamma_{*}$. 
We immediately get that the problem 
(\ref{equation})-(\ref{cond03}) has no solutions for $\gamma>\gamma_{*}.$ Indeed, in this case
the phase curve $C_{\gamma,\alpha}$ crosses the $v$-axis meaning that $f'$ vanishes and $f$
cannot be a solution.

Let us show that if $\alpha=\gamma_{*}^2u_{*}$ then $f$ is a bounded solution of 
(\ref{equation})-(\ref{cond03}).
As $C_{\gamma,\alpha}$ tends to the point $O$ as $s\rightarrow \infty $ tangentially with the 
line $L$, we have that for $t$ large enough $f^\prime(t)>0$ and $f^{\prime \prime}(t)<0$
which implies that $f$ is defined on the whole interval $[0,\infty)$. Furthermore
\begin{equation}
\frac{f^\prime(t)}{f(t)^2}\rightarrow 0\quad \text{and}\quad 
\frac{f^{\prime \prime}(t)}{f(t)f^\prime(t)}\rightarrow -(m+2)
\quad \text{as}\quad  t\rightarrow \infty . \label{tangentL}
\end{equation}
Hence $f^\prime(t)\rightarrow l\geqslant 0$ as $t\rightarrow \infty $ and
if we suppose $l>0$, from (\ref{tangentL}) we have
\[
f^{\prime \prime}(t)\sim -(m+2)l^2t \quad \text{as}\quad  t\rightarrow \infty ,
\]
and a contradiction with the fact that  
$f^\prime(t)\rightarrow l> 0$ as $t\rightarrow \infty $. So $l=0$
and $f$ is a solution to (\ref{equation})-(\ref{cond03}). Suppose now that $f$ 
is unbounded, i.e. $f(t)\rightarrow \infty $ as $t\rightarrow \infty $. Due to 
(\ref{tangentL}), there exists $t_{0}>0$ such that
\[
\forall t\geqslant t_{0}, \quad f^{\prime \prime}(t)\leqslant 
-\frac{m+2}{2}f(t)f^\prime(t).
\]
Integrating an dividing by $f(t)^2$ leads to
\[
\forall t \geqslant t_{0}, \quad 
\frac{f^\prime(t)}{f(t)^2}-\frac{f^\prime(t_{0})}{f(t)^2}
\leqslant -\frac{m+2}{4} \left ( 1-\frac{f(t_{0})^2}{f(t)^2} \right ).
\]
And using (\ref{tangentL}) leads to a contradiction as $t\rightarrow \infty $. 

Let us now look at what happens for $u_{-}\gamma^2<\alpha<u_{+}\gamma^2$. 
Because of the behavior of the vector field in the area 
$\left \{u>0\right \}\cap\left \{v<0\right \}$, we know that the phase
curve $C_{\gamma,\alpha}$ has to go to the singular point $O$ as $s\rightarrow \infty $ 
tangentially with the $u$-axis and below it. Thus, for large $t$ we have $f^\prime(t)>0$,
$f^{\prime \prime}(t)<0$ and again $f$ is defined on $[0,\infty).$ Moreover
\begin{equation}
\frac{f^\prime(t)}{f(t)^2}\rightarrow 0 \quad \text{and}\quad 
\frac{f^{\prime \prime}(t)}{f(t)f^\prime(t)}\rightarrow 0
\quad \text{as}\quad  t\rightarrow \infty . \label{tangentu}
\end{equation}
Hence $f^\prime(t)\rightarrow l\geqslant 0$ as $t\rightarrow \infty $ and supposing
$l>0$, we get from the following equality
$$f^{\prime \prime}(t)+(m+2)f(t)f^\prime(t)=-1-(m+2)\gamma \alpha
+3(m+1)\int_{0}^t f^\prime(\xi)^2 d \xi,$$
that
\[
f^{\prime \prime}(t)\sim (2m+1)l^2t \quad \text{as}\quad  t\rightarrow \infty ,
\]
which contradicts the fact that $f^\prime(t)\rightarrow l> 0$ as $t\rightarrow \infty $
and $f$ is a solution of (\ref{equation})-(\ref{cond03}).

Let us now prove that these solutions are unbounded. 
If $f$ were bounded, i.e. $f(t)\rightarrow \lambda $ as $t\rightarrow \infty $,
we can write (\ref{int_03}) with $\rho=t$ and $r=\infty$ in order to have
\[
-f^{\prime \prime}(t)f(t)+\frac{1}{2}f^\prime(t)^2
-(m+2)f^\prime(t)f(t)^2=(4m+5)\int_{t}^\infty f(\xi)f^\prime(\xi)^2d\xi. 
\]
Dividing this equality by $f^\prime(t)f(t)^2$ and using (\ref{tangentu}) leads to
\[
\frac{4m+5}{f^\prime(t)f(t)^2}\int_{t}^\infty f(\xi)f^\prime(\xi)^2d\xi \sim -(m+2)
\quad \text{as}\quad  t\rightarrow \infty,
\]
and a contradiction if $m\in \left [-\frac{5}{4},-1\right )$. 
If $m\in \left (-2,-\frac{5}{4} \right)$ we get 
\[
\int_{t}^\infty f(\xi)f^\prime(\xi)^2d\xi \sim -\frac{m+2}{4m+5}f^\prime(t)f(t)^2
\quad \text{as}\quad  t\rightarrow \infty,
\]
and using the fact that $f(t)\rightarrow \lambda$ as $t\rightarrow \infty $ we obtain
\[
\int_{t}^\infty f^\prime(\xi)^2d\xi \sim -\frac{m+2}{4m+5} \lambda f^\prime(t)
\quad \text{as}\quad  t\rightarrow \infty.
\]
We also have from (\ref{int_01})
\[
\int_{t}^\infty f^\prime(\xi)^2 d \xi=-\frac{1}{3(m+1)}\left (
f^{\prime \prime}(t)+(m+2)f(t)f^\prime(t) \right ).
\]
Combining these two equalities, we obtain
\[
\frac{f^{\prime \prime}(t)}{f(t)f^\prime(t)} \sim
-\frac{(m+2)^2}{4m+5}\not\neq 0,
\] 
and a contradiction with ($\ref{tangentu}$). We conclude that $f$ is an unbounded 
solution of (\ref{equation})-(\ref{cond03}).
\end{proof}

\begin{remark} For $-2<m<-1$ the critical value $\gamma_{*}$ depends on 
$m$, moreover $\gamma_{*}$ increases from $-\infty$ to $0$.
\end{remark}

\begin{remark}\label{equiv}
If $f$ is an unbounded solution of  $(\ref{equation})$-$(\ref{cond03})$ we can show that there
exists a positive constant $c$ such that
$$f(t)Ê\sim  ct^{\frac{m+2}{1-m}}Ê\quad \text{as}\quad t\rightarrow \infty. $$
For more details see {\rm \cite{BBequiv}}.
\end{remark}

\subsection{The case $m=-1$}
For $m=-1$, equation (\ref{equation}) reduces to 
\[%
\begin{array}
[c]{cl}
& f^{\prime\prime\prime}+ff^{\prime\prime}+f^{\prime2}=0\\
\Leftrightarrow & f^{\prime\prime\prime}+\left(  ff^{\prime}\right)  ^{\prime
}=0,
\end{array}
\]
and integrating on $[0,t]$ leads to
\begin{equation}
f''(t)+f(t)f'(t)=-1-\gamma f'(0). \label{ric1}
\end{equation}

Integrating $(\ref{ric1})$ and taking into account the boundary conditions (\ref{cond01})-(\ref{cond03})
leads to the Riccati equation%
\begin{equation}
f^{\prime}(t)+\dfrac{1}{2}f(t)^{2}=ct+d \label{Riccati}%
\end{equation}
with $c=-1-\gamma f^{\prime}\left(  0\right)  $ and $d=f^{\prime}(0)+\frac{\gamma^{2}}{2}.$

\begin{proposition}
For $m=-1$, solutions of $(\ref{equation})$-$(\ref{cond03})$ only exists if $\gamma<0.$ 
Moreover, if it is the case we have $f^{\prime}\left(  0\right)  \geq-\frac{1}{\gamma}.$
\end{proposition}
\begin{proof}
Suppose that $f$ is a solution of $(\ref{equation})$-$(\ref{cond03})$.
As $m<-\frac{1}{2}$ using proposition \ref{prop<-1/2} shows that $f^{\prime}(t)>0$
for all $t.$ Thus%
\[
\forall t\geq 0, \quad f^{\prime}(t)+\dfrac{1}{2}f(t)^{2}=ct+d\geq0
\]
and $c\geq0\Leftrightarrow-1\geq\gamma f^{\prime}(0),$ and we
deduce that $\gamma<0$ and $f^{\prime}\left(  0\right)  \geq-\frac{1}{\gamma}.$
\end{proof}

\begin{theorem}
Let $m=-1$, for every $\gamma<0$ the problem $(\ref{equation})$-$(\ref{cond03})$ admits an 
unique bounded solution with $f'(0)=-\frac{1}{\gamma}$ and an infinite number of 
unbounded solutions with $f'(0)>-\frac{1}{\gamma}$.
\end{theorem}
\begin{proof}
Let $f$ be the solution of (\ref{ivp}) with $\alpha\geq -\frac{1}{\gamma}$. From proposition 
\ref{prop<-1/2} we have that $f''<0$ and 
using equation (\ref{ric1}) we deduce that $f'$ cannot vanish. This implies that $f$ is defined on the whole 
interval $[0,\infty)$ and that $f'(t)$ has a limit $l\geq 0$ as $t \rightarrow \infty $. If we suppose 
$l>0$ we have that $f(t)f'(t)\rightarrow \infty $ as $t \rightarrow \infty $ and, using (\ref{ric1}), 
that $f''(t)\rightarrow -\infty $. Then $f'$ must become negative and this is a contradiction. Therefore 
$l=0$ and $f$ is a solution of (\ref{equation})-(\ref{cond03}). 

Suppose now that $f$ is bounded and writing (\ref{int_01}) with $\rho=0$ and $r=\infty$
we obtain that $f'(0)=-\frac{1}{\gamma}$ and the uniqueness.
\end{proof}

\begin{remark}
Let $f^{\prime}(0)=-\frac{1}{\gamma}$, then $(\ref{Riccati})$ can be integrated and we get that
the unique bounded solution of $\left( \ref{equation}\right)$-$(\ref{cond03})$ is given by%
\[
f(t)=\dfrac{2\sqrt{2d}}{\dfrac{\gamma-\sqrt{2d}}{\gamma
+\sqrt{2d}}e^{\sqrt{2d}t}-1}+\sqrt{2d},
\]
with $d=-\dfrac{1}{\gamma}+\dfrac{1}{2}\gamma^{2}$ and $\gamma<0.$
\end{remark}

\begin{remark}
Let $f$ be an unbounded solution of $\left( \ref{equation}\right)$-$(\ref{cond03})$.
Using $(\ref{int_01})$ with $m=-1$, $\rho=0$ and $r=t$ we obtain
$$f(t)f'(t)\rightarrow -(1+\gamma f'(0)) \quad \text{as}\quad  t\rightarrow \infty $$
from which we deduce that
$$f(t)\sim \sqrt{-2(1+\gamma f'(0))}\sqrt{t} \quad \text{as}\quad  t\rightarrow \infty.$$
\end{remark}
\subsection{The case $-1< m\leq-\frac{1}{2}$}
Let us introduce the following boundary value problem studied in \cite{brighi02}
\begin{equation}%
\left \{
\begin{array}
[c]{rl}%
&\hat g^{\prime\prime\prime}+\frac{n+1}{2} \hat g \hat g^{\prime\prime}-
n  \hat g^{\prime2}=0,\\
&\hat g(0)=0,\\
&\hat g^{\prime}\left(  0\right) =1,\\
&\hat g^{\prime}(\infty) =0.
\end{array}
\right .
\label{BB}%
\end{equation}

\begin{lemma}
For $n\in \left (-\frac{1}{3},\infty \right )$, the problem $(\ref{BB})$ admits as solution 
$\hat g$ which is increasing, strictly concave and verifies
$$\forall t\geq 0, \quad 0\leq \hat g(t)\leq \frac{2}{\sqrt{n+1}}.$$
\end{lemma}
\begin{proof}
See \cite{brighi02}.
\end{proof}

\begin{lemma}\label{g0}
For every $m\in (-1,\infty)$, there exists a function $g$ strictly concave and increasing 
that is a solution of the problem $(\ref{equation})$-$(\ref{cond03})$ with $\gamma=0$.
Moreover we have that
$$\forall t\geq 0, \quad 0\leq g(t)\leq \sqrt{\frac{2g'(0)}{m+2}}.$$
\end{lemma}
\begin{proof}
Let $m\in (-1,\infty)$ and let $\hat g$ a solution of the problem (\ref{BB}) 
with $n=\frac{2m+1}{3}$, then the function $g$ defined by 
$$g(t)=a.\hat g(bt)$$ 
with
$$a=\sqrt{\frac{g'(0)}3}\quad \text{and}\quad b=\sqrt{3g'(0)}$$
is a solution of (\ref{equation})-(\ref{cond03}).
\end{proof}

\begin{lemma}
For every $m\in (-1,\infty)$, the solutions $g$ of the problem $(\ref{equation})$-$(\ref{cond03})$ 
with $\gamma=0$ are such that
$$g'(0)\geq \frac{1}{\sqrt[3]{6(m+1)}}.$$
\end{lemma}
\begin{proof}
Let $m>-1$, using equality (\ref{int_01}) leads to
\[ g''(t)+(m+2)g(t)g'(t)+1=3(m+1)\int_{0}^{t}g'(s)^{2}ds,\]
and as $0<g'(t)\leq\alpha$ with $\alpha=g'(0)$ and $g\geq0$ we have
\[g''(t)+1\leq3(m+1)\alpha^{2}t.\]
Integrating this inequality we obtain
\[\forall t>0,\quad \frac{3(m+1)}{2}\alpha^{2}t^{2}-t+\alpha\geq0\]
and $\alpha\geq \frac{1}{\sqrt[3]{6(m+1)}}.$
\end{proof}
\bigskip

\begin{theorem} \label{ex-bounded}
Let $\gamma \in \mathbb{R}$. If $-1<m\leq -\frac{1}{2}$, the problem 
$(\ref{equation})$-$(\ref{cond03})$ admits a 
bounded solution $f$. This solution is positive at infinity, increasing, strictly concave and satisfies
$$\forall t\geq 0, \quad -\gamma \leq f(t) \leq \sqrt{\gamma^2+2\frac{f'(0)}{m+2}}.$$
Moreover if $\gamma \leq 0$ such a solution is unique.
\end{theorem}
\begin{proof}[Proof of existence]
Let $g$ be the solution of the problem (\ref{equation})-(\ref{cond03}) with $\gamma=0$
constructed in lemma \ref{g0}. 
\begin{itemize}
\item Case 1: $\gamma<0$. Since for all $k>0$ and all $t_{0}$ the function
\begin{equation}
f:\, t\rightarrow kg(kt+t_{0})\label{def_f}
\end{equation}
verifies (\ref{equation}) we want to choose $k$ and $t_{0}$ in order to have a
solution of (\ref{equation})-(\ref{cond03}) with $\gamma<0.$ First let us define the function
$h$ by
\begin{equation}
h:\, t\rightarrow\frac{g(t)^{3}}{g''(t)}. \label{fh}
\end{equation}
This function is well defined on $[0,\infty)$ and verify $h(0)=0$ and $h(t)\rightarrow -\infty$ 
as $t \rightarrow \infty$ because of proposition \ref{f_seconde} and lemma \ref{g0}.
Thus, there exists $t_{0}$ such
that $h(t_{0})=\gamma^3$ and choosing
$$k=-\frac{\gamma}{g(t_{0})},$$
wee see that for $\gamma<0$ the function $f$ defined by (\ref{def_f}) with these $k$ and $t_{0}$ is
a solution of (\ref{equation})-(\ref{cond03}).

\item Case 2:  $\gamma>0$. Let us consider again the function $h$ defined by (\ref{fh}). 
To use the previous 
method  we now have to look at $g(t)$ for the negative values of $t$. Let $(-T,\infty)$ be the 
maximal interval of existence of $g$. It is easy to see that if $g''$ does not vanish then $T=\infty$ 
because if $T<\infty$, in view of proposition \ref{blowT}, 
we have that $g(t)\rightarrow -\infty$, $g'(t)\rightarrow \infty$ and 
$g''(t)\rightarrow -\infty$ as $t \rightarrow -T$. Then as $m\leq -\frac{1}{2}$ equation 
(\ref{equation}) give $g'''(t)\rightarrow -\infty$, a contradiction. 

If $g''$ vanishes, let $t_{1}<0$ be such that $g''(t_{1})=0$ and
$g''<0$ on $(t_{1},0)$. Then $h$ is defined on $(t_{1},0]$
and $h(t)\rightarrow\infty$ as $t\rightarrow t_{1}$.

Suppose now that $g''<0.$ If $h$ is bounded on $(-\infty,0),$ there
exists $c>0$ such that
\[
\forall t<0,\quad \frac{g''(t)}{g(t)^{3}}>c.\]
Multiplying by $g(t)^{3}g'(t)$ and taking into account that $g<0$
leads to\[
g''(t)g'(t)<cg(t)^{3}g'(t).\]
 and integrating gives\[
\forall r<t<0,\quad g'(t)^{2}-g'(r)^{2}<\frac{c}{2}\left(g(t)^{4}-g(r)^{4}\right)\]
and finally\[
\forall r<t<0,\quad \frac{g'(t)^{2}}{g(r)^{4}}-\frac{g'(r)^{2}}{g(r)^{4}}<\frac{c}{2}
\left(\frac{g(t)^{4}}{g(r)^{4}}-1\right).\]
Since $g(r)\rightarrow -\infty $ as $r \rightarrow -\infty $, there exists $r_{0}<0$ such that
\[\forall r<r_{0},\quad \frac{g'(r)^{2}}{g(r)^{4}}>\frac{c}{4}\Leftrightarrow
\frac{g'(r)}{g(r)^{2}}>\frac{\sqrt{c}}{2}.\]
Integrating the last expression for $r<t<r_{0}$ we get\[
\forall r<t<r_{0},\quad -\frac{1}{g(t)}+\frac{1}{g(r)}>\frac{\sqrt{c}}{2}(t-r),\]
and passing to the limit as $r \rightarrow-\infty$ leads to a contradiction. Thus, $h$ is
unbounded on $(-\infty,0)$.

Therefore $h$ is always unbounded and there exists $t_{0}<0$ such
that $h(t_{0})=\gamma^{3}.$ If we choose $k=-\frac{\gamma}{g(t_{0})}$ the
function $f$ given by (\ref{def_f}) is a solution of (\ref{equation})-(\ref{cond03}).

From the boundedness of $g$ we deduce that $f$ is bounded too. Let $\lambda$ be the limit of $f$ at 
infinity. Using the boundedness and concavity of $f$ for large $t$ leads to
$$\underset{t\rightarrow \infty}{\lim}\, tf'(t)f(t)=
\underset{t\rightarrow \infty}{\lim}\,tf''(t)=0.$$
Writing (\ref{int_02}) with $\rho=0$ and $r=\infty$ leads to
\begin{equation}
f'(0)-\frac{m+2}{2}\left(\lambda^2-\gamma^2\right)=3(m+1)\int_{0}^{\infty}
\xi f'(\xi)^2d\xi>0. \label{finfini}
\end{equation}
%\newpage
And the result follows from the fact that $f$ is increasing. 
\end{itemize}
\vspace{-0.3cm}
\end{proof}

\begin{proof} [Proof of uniqueness]
Let $\gamma \leq 0$.
First let us remark that as $m\leq-\frac{1}{2}$ if $f$ is a solution of (\ref{equation})-(\ref{cond03}), 
$f$ is increasing and strictly concave. Thus we can define a function $v=v(y)$ such that
$$\forall t\geq 0,\quad v(f(t))=f'(t).$$
If $f$ is bounded, there exists $\lambda$ such that $f(t) \to \lambda$ as $t \to \infty$. Then $v$ is defined on $[-\gamma,\lambda)$, is positive and we have
\begin{align*}
f''(t)&=v(f(t))v'(f(t)), \\
f'''(t)&=v(f(t))v'(f(t))^2+v(f(t))^2v''(f(t)),
\end{align*}
and (\ref{equation}) leads to
\begin{equation}
\forall y\in [-\gamma,\lambda),\quad v''=-\frac{1}{v}\left (v'+(m+2)y\right)v'+(2m+1). 
\label{eq-v}
\end{equation}
We also have
$$v(-\gamma)=v(f(0))=f'(0)=\alpha>0,\quad v(\lambda):=
\underset{y\rightarrow \lambda}{\lim}\, v(y)
=\underset{t\rightarrow \infty}{\lim}\, f'(t)=0  \quad \text{and}\quad 
v'(-\gamma)=-\frac{1}{\alpha}.$$
Suppose now that there are two bounded solutions $f_{1}$ and $f_{2}$ of 
(\ref{equation})-(\ref{cond03}) and let $\lambda_{i}$ be the limit of $f_{i}$ at infinity for $i=1,2.$ They give $v_{1},\, v_{2}$ solutions of equation (\ref{eq-v}) defined respectively on 
$[-\gamma,\lambda_{1})$ and $[-\gamma,\lambda_{2})$ such that
$$v_{1}(-\gamma)=\alpha_{1}, \, v_{2}(-\gamma)=\alpha_{2}, \quad \text{and}\quad  
v_{1}(\lambda_{1})=v_{2}(\lambda_{2})=0.$$
Let us suppose that $\alpha_{1}<\alpha_{2}$ and show that $\lambda_{1}\leq \lambda_{2}$. 
If, on the contrary, $\lambda_{1}>\lambda_{2}$ the function $w=v_{1}-v_{2}$ verifies 
$w(-\gamma)<0$, $w(\lambda_{2})=v_{1}(\lambda_{2})>0$ and 
$w'(-\gamma)=\frac{\alpha_{1}-\alpha_{2}}{\alpha_{1}\alpha_{2}}<0$. Then $w$ admits 
a negative minimum at some point $x\in(-\gamma,\lambda_{2})$. So we have 
$v_{1}(x)<v_{2}(x)$, $v'_{1}(x)=v'_{2}(x)$ and $v''_{1}(x)\geq v''_{2}(x)$. We also have
\begin{equation}
v''_{1}(x)-v''_{2}(x)=\left (\frac{1}{v_{2}(x)}-\frac{1}{v_{1}(x)} \right )
\left (v'_{1}(x)+(m+2)x\right)v'_{1}(x), 
\label{eq-v1}
\end{equation}
and
\begin{equation}
\left (v'_{1}(x)+(m+2)x\right)v'_{1}(x)=
\left (f''_{1}(t)+(m+2)f_{1}(t)f'_{1}(t)\right)\frac{f''_{1}(t)}{f'_{1}(t)^2}, 
\label{eq-v2}
\end{equation}
with $t$ such that $x=f_{1}(t)$. As $f_{1}$ is bounded, writing (\ref{int_01}) with $\rho=t$
and $r=\infty$ leads to
\begin{equation}
f''_{1}(t)+(m+2)f_{1}(t)f'_{1}(t)=-3(m+1)\int_{t}^{\infty}f'_{1}(\xi)^2d\xi<0. \label{eq-v3}
\end{equation}
Using this inequality and the fact that $f''_{1}(t)<0$, (\ref{eq-v1}) and (\ref{eq-v2}) leads to 
$v''_{1}(x)<v''_{2}(x)$ and a contradiction. Therefore we have $\lambda_{1}\leq \lambda_{2}$.

Now let us prove that $v_{1}\leq v_{2}$ on 
$[-\gamma,\lambda_{1}).$ For that suppose there exists some $y\in (-\gamma,\lambda_{1})$ such 
that $v_{1}(y)>v_{2}(y)$ and set $w=v_{1}-v_{2}$. As $\alpha_{1}<\alpha_{2}$, 
$w(-\gamma)<0$ and using the fact that $w(\lambda_{1})\leq 0$ we deduce that $w$ admits a 
positive maximum at a point $x\in (-\gamma,\lambda_{1})$. Thus  $v_{1}(x)>v_{2}(x)$, 
$v'_{1}(x)=v'_{2}(x)$ and $v''_{1}(x)\leq v''_{2}(x)$. 

Using inequality (\ref{eq-v3}) and the fact that $f''_{1}(t)<0$, (\ref{eq-v1})-(\ref{eq-v2}) leads to 
$v''_{1}(x)>v''_{2}(x)$ and a contradiction. Therefore we have $v_{1}\leq v_{2}$ on 
$[-\gamma,\lambda_{1})$ and
$$\int_{0}^{\infty}f'_{1}(\xi)^2d\xi=\int_{-\gamma}^{\lambda_{1}}v_{1}(y)dy
<\int_{-\gamma}^{\lambda_{1}}v_{2}(y)dy
\leq \int_{-\gamma}^{\lambda_{2}}v_{2}(y)dy
=\int_{0}^{\infty}f'_{2}(\xi)^2d\xi.$$
Since
$$-1-(m+2)\gamma\alpha_{i}=-3(m+1)\int_{0}^{\infty}f'_{i}(\xi)^2d\xi$$
we get 
$$\gamma(\alpha_{1}-\alpha_{2})<0$$
and as $\alpha_{1}-\alpha_{2}<0$ this leads to $\gamma>0$ and a contradiction.
\end{proof}

\begin{remark}
Let  $m\in \left (-1,-\frac{1}{2} \right ]$ and $f$ a bounded solution of $(\ref{equation})$-$(\ref{cond03})$.
As $f$ is strictly concave on $[0,\infty)$ we have
$rf'(r)<f(r)+\gamma$ for $r>0$. If $\lambda$ denotes the limit of $f$ at infinity we get
$$\int_{0}^{\infty}rf'(r)^2dr<\frac{1}{2}(\lambda+\gamma)^2.$$
Then as $f'(0)>0$, $(\ref{finfini})$ becomes
$$ (4m+5)\lambda^2+6(m+1)\gamma \lambda+(2m+1)\gamma^2>0$$
and for $\gamma>0$ we have
$$-\frac{2m+1}{4m+5}\gamma<\lambda \leq \sqrt{\gamma^2+2\frac{f'(0)}{m+2}}.$$
\end{remark}
\bigskip

\begin{theorem}
Let $-1< m\leq-\frac{1}{2}$, then for $\gamma<0$ the 
problem $(\ref{equation})$-$(\ref{cond03})$ admits many infinitely unbounded solutions.
\end{theorem}
\begin{proof}
We follow an idea of \cite{guedda}. Consider the initial value problem (\ref{ivp}) 
with $\gamma<0$ and let $f_{\alpha}$ be its solution on 
$[0,T_{\alpha})$. Writing (\ref{int_01}) with $\rho=0$ and $r=t<T_{\alpha}$ leads to
\begin{equation}
f''_{\alpha}(t)+(m+2)f_{\alpha}(t)f'_{\alpha}(t)
=-(m+2)\gamma\alpha-1+3(m+1)\int_{0}^tf'_{\alpha}(\xi)^2d\xi. \label{intf2}
\end{equation}
For the remainder of the proof let us choose $\alpha\geq -\frac{1}{(m+2)\gamma}$. Then
\begin{equation}
\forall t \in [0,T_{\alpha}), \quad  f''_{\alpha}(t)+(m+2)f_{\alpha}(t)f'_{\alpha}(t)>0,
\label{vanish}
\end{equation}
and it follows that $f'_{\alpha}(t)>0$ for all $t$ in $[0,T_{\alpha})$. Indeed, since 
$f'_{\alpha}(0)=\alpha>0$, we 
should have $f''_{\alpha}(t_{1})\leq 0$ for $t_{1}$ the first point where $f'_{\alpha}(t)$
vanishes that leads to a contradiction with (\ref{vanish}). Using lemma
\ref{concavite} we have that $f''_{\alpha}<0$ on $[0,T_{\alpha})$, then 
$f'_{\alpha}(t)\rightarrow l\in [0,\infty)$ as $t\rightarrow T_{\alpha}$. As $f_{\alpha}$ is 
strictly concave and increasing we deduce that $T_{\alpha}=\infty$ and $f_{\alpha}(t)>0$ on
$[0,\infty).$

If $l\not\neq 0$ we have $f_{\alpha}(t)\sim lt$ as $t\rightarrow \infty,$
and using (\ref{intf2}) we obtain that $f''_{\alpha}(t) \sim -(4m+5)l^2t$ as $t\rightarrow \infty$
that is a contradiction with $f'_{\alpha}(t)\sim l$ as $t\rightarrow \infty.$

Finally we get $l=0$ and $f_{\alpha}$ verifies (\ref{equation})-(\ref{cond03}). Furthermore,
from (\ref{intf2}) and the choice of $\alpha$ we deduce that $f_{\alpha}$ is unbounded.
\end{proof}

\begin{remark} As in the case $-2<m<-1$, we also have that if $f$ is an unbounded 
solution of $(\ref{equation})$-$(\ref{cond03})$, there exists a positive constant $c$
such that 
$$f(t) \sim ct^{\frac{m+2}{1-m}}\quad \text{as} \quad t\rightarrow \infty .$$
For more details see {\rm \cite{BBequiv}}.
\end{remark}

\begin{remark}
For $m=-\frac{1}{2}$, equation $(\ref{equation})$ reduces to
$$f^{\prime\prime\prime}+\dfrac{3}{2}ff^{\prime\prime}=0$$
which is the Blasius equation. This equation is investigated in 
{\rm $\cite{Cop}$} and {\rm $\cite{Hart}$} 
and its concave solutions are studied in {\rm $\cite{brighi03}$} and {\rm $\cite{Ish}.$} 
See also {\rm $\cite{new}$}.
\end{remark}

\subsection{The case $m>-\frac{1}{2}$}

\begin{theorem} \label{thbb1} Let  $\gamma \in \mathbb{R}$. For any $m\geq -\frac{1}{2}$ the problem 
$(\ref{equation})$-$(\ref{cond03})$ admits one and only one concave solution $f$ which is positive 
at infinity and such that
\begin{equation}
\forall t \geq 0, \quad -\gamma \leq f(t) \leq \sqrt{\gamma^2+2\frac{f'(0)}{m+2}}.
\end{equation}
\end{theorem}
\begin{proof}[Proof of existence]
Let $g$ be the solution of (\ref{equation})-(\ref{cond03}) with $\gamma=0$ constructed in lemma 
\ref{g0}.
\begin{itemize}
\item Case 1: $\gamma<0$. The same proof as in the theorem \ref{ex-bounded} works 
well in this case too.
\item Case 2: $\gamma>0$. As in theorem \ref{ex-bounded} we denote by $(-T,\infty)$ the 
maximal interval of existence of $g$ and we again consider the function $h$ defined by (\ref{fh}). 
Using lemma \ref{concavite}, 
$g$ is strictly concave, increasing and $h$ is defined on $(-T,\infty)$. Let us prove that $h$ is 
unbounded on $(-T,\infty)$.

If $T=\infty$ the reasoning used for theorem \ref{ex-bounded} still works, so let us suppose that $T<\infty$.
Using proposition \ref{blowT} we have that $g(t)\rightarrow -\infty$, $g'(t)\rightarrow \infty$ and 
$g(t)''\rightarrow -\infty$ as $t \rightarrow -T$. Differentiating (\ref{equation}) leads to
\begin{equation}
\left (g'''e^{(m+2)G} \right )'=3me^{(m+2)G}g'g'' \label{g3}
\end{equation}
with $G$ any anti-derivative of $g$. Then, as $g'''(0)=(2m+1)g'(0)^2$ using (\ref{g3}) we have 
that $g'''>0$ on $(-T,\infty)$ and setting $\beta=\frac{2m+1}{m+2}$ leads to
$$-gg''+\beta g'^2>0.$$
We deduce that the function $\phi=g'(-g)^{-\beta}$ is positive and increasing on $(-T,0)$ and that
$\phi$ is bounded as $t\rightarrow -T$. If $h$ is bounded on $(-T,0)$, there exists a 
positive constant $c$ such that $h(t)^{-1}>c>0$ and we have that
$$\forall t< 0, \quad g''(t)g'(t)<g(t)^3g'(t).$$
Integrating leads to
$$\forall r<t<0, \quad -g'(r)^2<g'(t)^2-g'(r)^2<\frac{c}{2}\left(g(t)^4-g(r)^4 \right)$$
and
$$\forall r<t<0, \quad -\frac{g'(r)^2}{g(r)^4}<\frac{c}{2}\left(\frac{g(t)^4}{g(r)^4} -1\right).$$
If we let $t$ going to zero we obtain that 
$$\frac{g'(r)}{g(r)^2}\geq \sqrt{\frac{c}{2}}$$
and
$$0<\sqrt{\frac{c}{2}}\leq \phi(r)(-g(r))^{\beta-2}\rightarrow 0 \quad \text{as}
\quad r\rightarrow -T$$
because $\beta<2$. This is a contradiction. 
\end{itemize}
As in any case $h$ is unbounded we conclude the same way as in theorem \ref{ex-bounded}.
%\vspace{-0.3cm}
\end{proof}

%\newpage
\begin{proof}[Proof of uniqueness]
Let $f_{1}$ and $f_{2}$ be two concave solutions of (\ref{equation})-(\ref{cond03}) such that 
$f'_{1}(0)>f'_{2}(0)$ and let $k=f_{1}-f_{2}$. The function $k$ verify $k(0)=0$, $k'(0)>0$, $k''(0)=0$
and $k'(\infty)=0$. Moreover, using proposition \ref{prop>-1/2} we have $f'_{1}(0)>0$, $f'_{2}(0)>0$ 
and  
$$k'''(0)=(2m+1)\left(f'_{1}(0)+f'_{2}(0)\right)\left(f'_{1}(0)-f'_{2}(0)\right)>0.$$
Then, the function $k$ is convex near $0$ and there exists $t_{0}>0$ such that $k'(t)>0$ on $(0,t_{0}]$, 
$k''(t_{0})=0$, $k'''(t_{0})\leq 0$ and $k(t_{0})>0$.

Using the fact that $f''_{1}(t_{0})=f''_{2}(t_{0})$ we obtain
$$k'''(t_{0})=(2m+1)k'(t_{0})\left(f'_{1}(t_{0})+f'_{2}(t_{0})\right)
-(m+2)f''_{1}(t_{0})k(t_{0})>0$$
wich leads to a contradiction with $k'''(t_{0})\leq 0$.
\end{proof}

\begin{lemma} \label{domain}
Let $m>-\frac{1}{2}$ and $f$ be a concave-convex solution of $\left(  \ref{equation}%
\right)  $-$\left(  \ref{cond03}\right)$. Let $t_{0}$ be the point such that
$f^{\prime\prime}(t_{0})=0$, then the curve $s\mapsto\left(  u(s),v(s)\right)
$ defined by $\left(  \ref{new_function}\right)  $ with $\tau=t_{0}$ is a
positive semi-trajectory which lies in the bounded domain%
\[
\mathcal{D}_{-}=\left\{  \left(  u,v\right)  \in\mathbb{R}^{2}\, ;\quad 
-\frac{m+2}{2}<u<0 \quad \text{and} \quad 0\leq v<-\left(  m+2\right)  u\right\}  .
\]
\end{lemma}
\begin{proof}
In view of proposition \ref{prop>-1/2} we know that $f$ is positive, decreasing and
convex on $\left[  t_{0},\infty\right)  ,$ thus%
\begin{equation}
\forall t\geq t_{0,}\quad \frac{f^{\prime}(t)}{f(t)^{2}}<0 \quad \text{and}\quad %
\frac{f^{\prime\prime}(t)}{f(t)^{3}}>0. \label{i3}%
\end{equation}
As $f$ is bounded, writing $\left(  \ref{int_01}\right)  $ with $\rho=t$
and $r=\infty$ we have
\begin{equation}
f^{\prime\prime}(t)+(m+2)f(t)f^{\prime}(t)=-3(m+1)\int_{t}^{\infty}f^{\prime
}(\xi)^{2}d\xi<0, \label{i1}%
\end{equation}
and if we denote by $\lambda$ the limit of $f$ at infinity, integrating leads
to%
\begin{equation}
f^{\prime}(t)+\frac{m+2}{2}f(t)^{2}>\frac{m+2}{2}\lambda^{2}\geq0. \label{i2}%
\end{equation}
\bigskip From $\left(  \ref{i1}\right)  $ and $\left(  \ref{i2}\right)  $ we
obtain that%
\begin{equation}
\frac{f^{\prime}(t)}{f(t)^{2}}+\frac{m+2}{2}>0 \quad \text{and} \quad \frac{f^{\prime
\prime}(t)}{f(t)^{3}}+(m+2)\frac{f^{\prime}(t)}{f(t)^{2}}<0, \label{i4}%
\end{equation}
\bigskip and%
\[
\forall t\geq t_{0,}\quad f(t)\geq\frac{1}{\frac{m+2}{2}\left(
t-t_{0}\right)  +\frac{1}{f(t_{0})}}%
\]
which implies%
\[
\int_{t_{0}}^{\infty}f(\xi)d\xi=\infty.
\]
Hence the trajectory $s\mapsto\left(  u(s),v(s)\right)  $ is defined on the
whole interval $\left[  0,\infty\right)  $ and using $\left(  \ref{i3}\right)
$ and $\left(  \ref{i4}\right)  $ leads to the result.
\end{proof}

\begin{remark}
For $m=1$, equation $(\ref{equation})$ reduces to
\[
f^{\prime\prime\prime}+3ff^{\prime\prime}-3f^{\prime2}=0.
\]
Let $f=g+\eta$ with $\eta>0$, we have%
\[
g^{\prime\prime\prime}+3\eta g^{\prime\prime}=3g^{\prime2}-3gg^{\prime
\prime}.
\]
Solving $g^{\prime\prime\prime}+3\eta g^{\prime\prime}=0$ with
$g(0)=-\gamma-\eta,$ $g^{\prime}(\infty)=0$ and $g^{\prime\prime}(0)=-1$
leads to%
\[
g(t)=-\dfrac{1}{9\eta^{2}}\left(  e^{-3\eta t}-1\right)  -\gamma-\eta
\]
and if we choose $\eta$ as the unique positive number such that
$9\eta^{3}+9\gamma\eta^{2}-1=0$ we easily see that $g$ satisfies
$g'^2-gg''=0.$
It follows that $f$ given by
\[
f(t)=-\dfrac{1}{9\eta^{2}}\left(  e^{-3\eta t}-1\right)  -\gamma
\]
is a solution of $(\ref{equation})$-$(\ref{cond03})$. Moreover, since 
$f^{\prime\prime}(t)=-e^{-3\eta t}<0,$ this is the unique concave solution of 
$(\ref{equation})$-$(\ref{cond03})$.
\end{remark}

\begin{theorem}
\label{t4}
Let $m\in \left (-\frac{1}{2},1\right ]$, then for any $\gamma \in \mathbb{R}$ the problem 
$(\ref{equation})$-$(\ref{cond03})$ admits one and only one solution, which is concave.
\end{theorem}
\begin{proof}
Taking into account proposition \ref{prop>-1/2} and theorem \ref{thbb1}, we just have to
consider the case $m\in \left (-\frac{1}{2},1 \right ]$ and prove that in this case concave-convex solutions cannot exist. 

\[
%TCIMACRO{\FRAME{itbpFU}{10.7989cm}{8.1166cm}{0cm}{\Qcb{$m<-2$}}{}%
%{Figure}{\special{ language "Scientific Word";  type "GRAPHIC";
%maintain-aspect-ratio TRUE;  display "USEDEF";  valid_file "T";
%width 10.7989cm;  height 8.1166cm;  depth 0cm;  original-width 7.6803in;
%original-height 5.76in;  cropleft "0";  croptop "1";  cropright "1";
%cropbottom "0";  tempfilename 'm_inf_-2.jpg';tempfile-properties "XPR";}}}%
%BeginExpansion
\raisebox{-0cm}{\parbox[b]{10.7989cm}{\begin{center}
\includegraphics[scale=.4]
{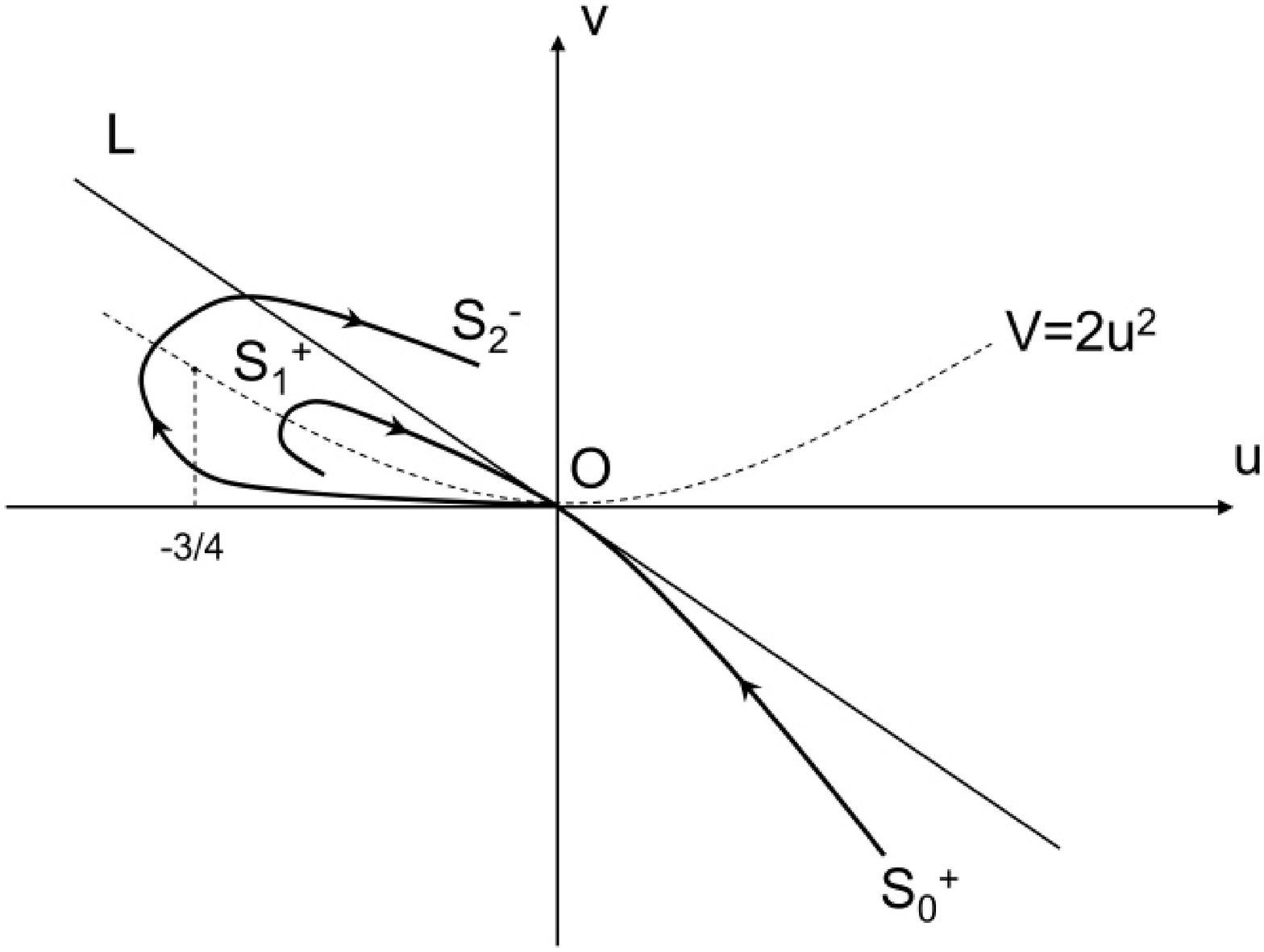}
\\
$-\frac{1}{2}<m<1$\\
{Fig 3.5.1}
\end{center}}}%
%EndExpansion
%\tag{($m<-2$)}%
\]

Suppose that $f$ is a concave-convex solution of (\ref{equation})-(\ref{cond03}) and denote 
by $t_0$ the point where $f^{\prime \prime}(t_0)=0$. Consider the positive semi-trajectory 
$s\longmapsto(u(s),v(s))$ defined in lemma \ref{domain}, we have
$$u(0)=\frac{f'(t_0)}{f(t_0)^2}<0 \quad \text{and} \quad v(0)=0.$$
Refering to Fig 2.2.1 we see that the behavior of the corresponding phase curve is related to the one of the separatrices $S_{2}^-$ and $S_{1}^+$.

As $s$ increases, the separatrix $S_{2}^-$ leaves the singular point $O$ to the left 
tangentially with $L_{0}$, and either does not cross the isocline $P(u,v)=0$, or crosses it through 
a point $(u_{2},2u_{2}^2)$ such that $u_{2}\leq  -\frac{3}{4}$ and next intersects the straight line $L$.

As $s$ decreases, the separatrix $S_{1}^+$ leaves the singular point $O$ 
to the left tangentially with $L$ and crosses the isocline $P(u,v)=0$ through a point 
$(u_{1},2u_{1}^2)$ such that $-\frac{3}{4} \leq u_{1}<0$ and next stays in the bounded region 
$\mathcal{D}_{-}$ (see Fig 3.5.1).

 In view of the behavior of the separatrices we see that this semi-trajectory cannot remain 
 in the  bounded domain ${\cal D}_-$ and a contradiction.
\end{proof} 

\begin{theorem} \label{concave-convex}
Let $m>1$, then for any $\gamma\in \mathbb{R}$ the problem $(\ref{equation})$-$(\ref{cond03})$ 
has infinitely many concave-convex solutions.
\end{theorem}
\begin{proof}
\begin{itemize}
\item Case 1: Let $\gamma<0$, and consider the initial value problem $\mathcal{P}_{m,\gamma,\alpha}$ 
given by ($\ref{ivp}$) and the corresponding phase curve $C_{\gamma,\alpha}$ of the system 
(\ref{system}) defined by (\ref{new_function}) with $\tau =0$. The separatrices we are concerned with, are $S_0^+$, $S_1^+$ and $S_2^-$.
\[
%TCIMACRO{\FRAME{itbpFU}{10.7989cm}{8.1166cm}{0cm}{\Qcb{$m<-2$}}{}%
%{Figure}{\special{ language "Scientific Word";  type "GRAPHIC";
%maintain-aspect-ratio TRUE;  display "USEDEF";  valid_file "T";
%width 10.7989cm;  height 8.1166cm;  depth 0cm;  original-width 7.6803in;
%original-height 5.76in;  cropleft "0";  croptop "1";  cropright "1";
%cropbottom "0";  tempfilename 'm_inf_-2.jpg';tempfile-properties "XPR";}}}%
%BeginExpansion
\raisebox{-0cm}{\parbox[b]{10.7989cm}{\begin{center}
\includegraphics[scale=.4]
{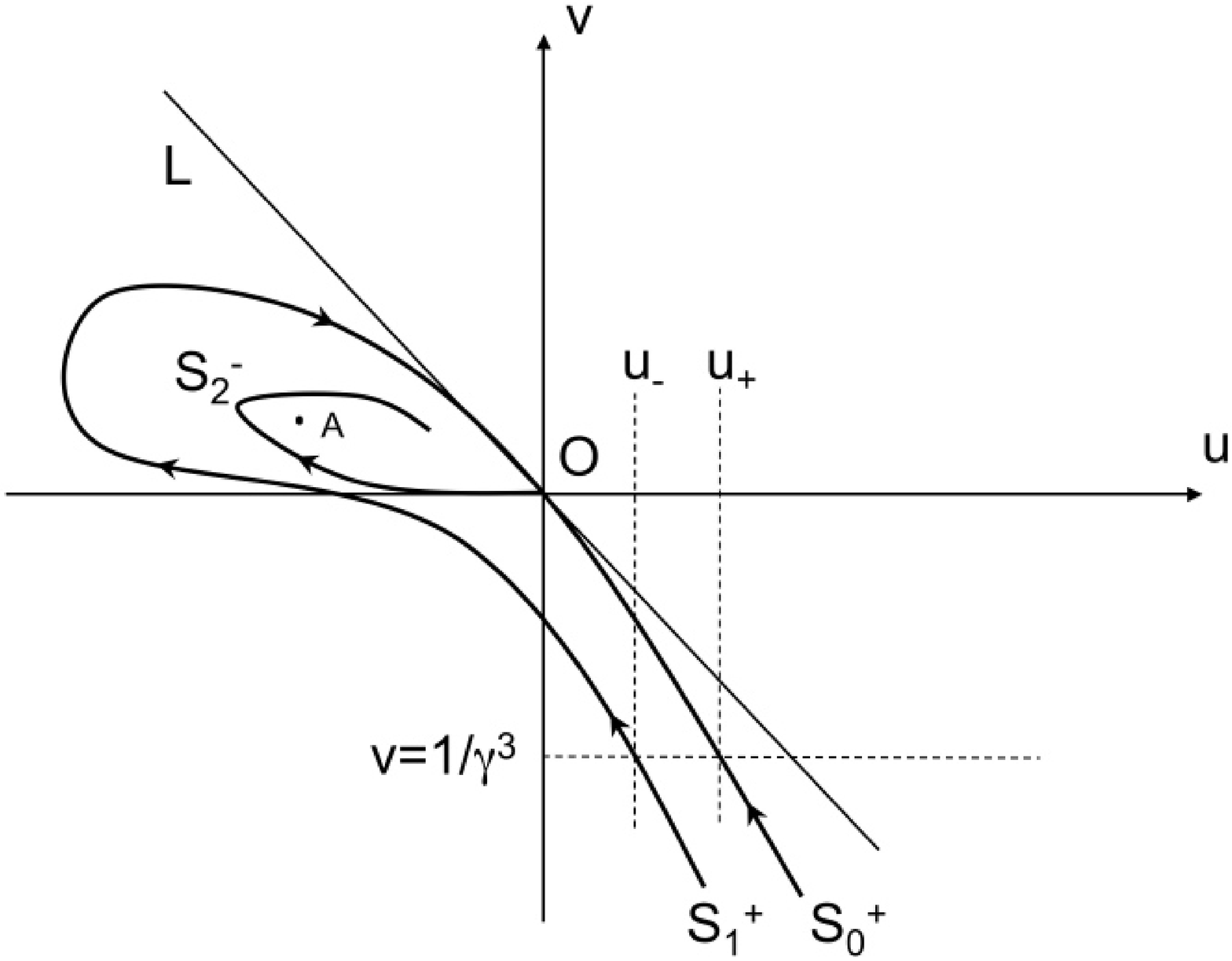}
\\
$m>1$\\
{Fig 3.5.2}
\end{center}}}%
%EndExpansion
%\tag{($m<-2$)}%
\]

As $s$ increases, the separatrix $S_{2}^-$ leaves the singular point $O$ to the left 
tangentially with $L_{0}$, and crosses the isocline $P(u,v)=0$ through a point $(u_{2},2u_{2}^2)$ 
such that $-\frac{3}{4}\leq u_{2}<0$ and then stay in the bounded region $\mathcal{D}_{-}$.

As $s$ decreases, the separatrix $S_{1}^+$ leaves the singular point $O$ 
to the left tangentially with $L$ and crosses the isocline $P(u,v)=0$ through a point 
$(u_{1},2u_{1}^2)$ such that $u_{1} \leq -\frac{3}{4}$. Then it intersects successively the $u$-axis 
and the $v$-axis and next stays in the region $\{u>0\}\cap \{v<0\}$ and goes to infinity with  a 
slope that stays between $-3u-(m+2)$ and $-(m+2)$.

As $s$ decreases, the separatrix $S_{0}^+$ leaves the singular point $O$ 
to the right tangentially with $L$ and below $L$. Then it stays in the region 
$\{u>0\}\cap \{v<0\}$ and goes to infinity (see Fig 3.5.2).

Looking at these separatrices we see that the straight line 
$v=\frac{1}{\gamma^3}$ crosses $S_{0}^+$ and $S_{1}^+$ through two points 
$(u_{-},\frac{1}{\gamma^3})$ and $(u_{+},\frac{1}{\gamma^3})$ with $u_{-}<u_{+}$.

For $\alpha\in[\gamma^2u_-,\gamma^2u_+)$, the trajectory $C_{\gamma,\alpha}$
intersects the $u$-axis for some $s_0$ and remains in the domain defined by the separatrix $S_1^+$
for $s>s_0$. It follows from the Poincar\'e-Bendixson Theorem that $C_{\gamma,\alpha}$ is a positive
semi-trajectory whose $\omega$-limit set is the point $O$ if $\alpha=\gamma^2u_-$, and either the 
singular point $A$ or a limit cycle surrounding $A$ if $\gamma^2u_-<\alpha<\gamma^2u_+$.
Therefore $f$ is positive as long it exists. 
Since $F_m(u,0)=-(m+\frac{1}{2})<0$ by (\ref{Fm}),
such a limit cycle cannot cross the $u$-axis and there exists $t_0>0$ such that $f'(t)<0$ 
and $f''(t)>0$ for $t>t_0$. Hence $f$ is defined on $[0,\infty)$, $f'(t)\to l\leq 0$ as $t\to\infty$ 
and if we suppose that $l<0$ 
we get a contradiction with the positivity of $f$. Consequently, if $\alpha\in
[\gamma^2u_-,\gamma^2u_+)$ then $f$ is a concave-convex solution of 
(\ref{equation})-(\ref{cond03}). To complete the proof in this case, let us remark that
for $\alpha\notin[\gamma^2u_-,\gamma^2u_+]$, in view of lemma \ref{domain}, the function $f$ cannot be a solution of (\ref{equation})-(\ref{cond03}), and that for 
$\alpha=\gamma^2 u_{+}$ $f$ is the concave solution.
\item Case 2: Let $\gamma\geq 0$ and $g$ be a concave-convex solution of 
$(\ref{equation})$-$(\ref{cond03})$ with $g(0)>0$ and $g'(0)>0$. 
Such a solution exists due to the precedent case. The function $g$ is defined on 
$(-T,\infty]$ and is strictly concave on $(-T,0]$ by lemma \ref{concavite}. Then, as 
$g'(0)>0$, there exists $t_{1}<0$
such that $g(t_{1})=0$. We know that for all $k>0$
and all $t_{0}$ the function $f(t)=kg(kt+t_{0})$ verifies $(\ref{equation})$ and we want to choose 
$k$ and $t_{0}$ to obtain a solution of $(\ref{equation})$-$(\ref{cond03})$ with $\gamma \geq 0$.

Let us consider again the function $h$ defined by $(\ref{fh})$. As  $g''$ does not vanish on 
$(-T,t_{1}]$
$h$ exists on $(-T,t_{1}]$, verifies $h(t_{1})=0$ and is unbounded. Indeed, to prove that $h$ is unbounded, 
we use the same proof as in theorem \ref{ex-bounded} if $T=\infty$ and 
the same as in theorem \ref{thbb1} if $T<\infty$. Then we construct a solution of 
$(\ref{equation})$-$(\ref{cond03})$
with $\gamma \geq 0$ by setting $k=-\frac{\gamma}{g(t_{0})}$ and the proof is complete.
\end{itemize}
\vspace{-0.3cm}
\end{proof}

\begin{remark}
Suppose given $\gamma<0$
\begin{itemize}
\item
As $u_{+}$ is the intersection of the separatrix $S_{0}^+$ that lies in the domain 
$\{u>0\} \cap \{v<0\}$ and the straight line $v=\frac{1}{\gamma^3}$ with 
$\gamma<0$, we have $u_{+}>0$.
\item
If $\gamma$ is such that $u_{-}>0$, then all the concave-convex solutions of the problem 
$(\ref{equation})$-$(\ref{cond03})$ are increasing-decreasing.
\item
If $\gamma$ is such that $u_{-}<0$, then for $\alpha\in[\gamma^2 u_{-},0]$ we get 
concave-convex solutions of $(\ref{equation})$-$(\ref{cond03})$ which are 
decreasing, and for $\alpha\in(0,\gamma^2 u_{+})$ we get concave-convex solutions 
increasing-decreasing.
\end{itemize}
\end{remark}

\begin{proposition}
Let $m>1$, then for every $\gamma \in \mathbb{R}$ there is an unique concave-convex solution 
that verify $f(t)\rightarrow l>0$ as $t\rightarrow \infty$ and all the other 
concave-convex solutions are such that $f(t)\rightarrow 0$ as $t\rightarrow \infty$.
\end{proposition}
\begin{proof}
Let $m>1$ and let $f$ be a concave-convex solution of $(\ref{equation})$-$(\ref{cond03})$.
Since $f$ is positive and decreasing at infinity, $f(t)\rightarrow\lambda\geq 0$ as $t\to\infty$.
If $f$ corresponds to the separatrix $S_1^+$ (i.e. $f'(0)=\gamma^2u_{-}$) then we prove as in 
proposition \ref{limf} that $\lambda>0$, and if $\gamma^2u_-<f'(0)<\gamma^2u_+$, there exists
$c>0$  such that $\vert f'(t)\vert>c\vert f(t)^2\vert$ for $t$ large enough, in such a way that
$\lambda=0$. 
\end{proof}

\begin{remark}
\label{r9}
For $1<m<\frac{3}{2}$  the singular point $A$ is an unstable focus,
which implies that at least one cycle surrounding $A$ has to exist. If $m>\frac{3}{2}$ then $A$
is attractive and it seems that cycles do not exist. If it is the case, we have
$$\frac{f'(t)}{f(t)^2}\sim-\frac{1}{2}\quad \text{and}\quad \frac{f''(t)}{f(t)^3}\sim
\frac{1}{2}\quad \text{as} \quad t\rightarrow \infty,$$
which easily give
$$f(t)\sim \frac{2}{t}\quad \text{as}\quad t\rightarrow \infty.$$
 \end{remark}

\section{Conclusion}

\bigskip

\begin{itemize}

\item For $m<-2$ there exists 
$\gamma_{*}>\sqrt[3]{\frac{2}{(m+2)^2}}$ such that the problem
$\left(  \ref{equation}\right)  $-$\left(  \ref{cond03}\right)  $ has
no solution for $\gamma<\gamma_{\ast}$, one and only one solution
for $\gamma=\gamma_{\ast}$ and infinitely many solutions for $\gamma>\gamma_{\ast}.$

For $\gamma=\gamma_{*}$ we have that $f(t)\rightarrow \lambda <0$
as $t\rightarrow \infty$ and for every $\gamma> \gamma_{*}$ there are two solutions $f$
such that $f(t)\rightarrow \lambda <0$ as $t\rightarrow \infty$ and all the other solutions verify 
$f(t)\rightarrow 0$ as $t\rightarrow \infty$

Moreover, if $f$ is a solution of $\left(  \ref{equation}\right)  $-$\left(
\ref{cond03}\right)  $, then $f$ is negative, strictly concave and
increasing.

\item For $m=-2$ and for every $\gamma\in\mathbb{R}$, the problem $\left(
\ref{equation}\right)  $-$\left(  \ref{cond03}\right)  $ has no solution.

\item For $-2<m<-1$, there exists $\gamma_{\ast}<0$ such that the problem
$\left(  \ref{equation}\right)  $-$\left(  \ref{cond03}\right)  $ has
no solution for $\gamma>\gamma_{\ast}$, one and only one solution which is bounded 
for $\gamma=\gamma_{\ast}$ and two bounded solutions and infinitely many 
unbounded solutions for $\gamma<\gamma_{\ast}.$

Moreover, if $f$ is a solution of $\left(  \ref{equation}\right)  $-$\left(
\ref{cond03}\right) $, then $f$ is positive, strictly concave,
increasing and $f^{\prime}(0)\geq -\frac{1}{(m+2)\gamma}.$

\item For $m=-1$ the problem $\left(  \ref{equation}\right)  $-$\left(
\ref{cond03}\right) $ only admits solutions for $\gamma<0$. In this case there is an unique 
bounded solution with $f'(0)=-\frac{1}{\gamma}$ and an infinite 
number of unbounded solutions with $f'(0)>-\frac{1}{\gamma}$. Moreover all the
solutions are positive, strictly concave and increasing.

\item For $-1< m<-\frac{1}{2}$ the problem $\left(  \ref{equation}\right)  $-$\left(
\ref{cond03}\right) $ admits at least one bounded solution for $\gamma \in \mathbb{R}$ and 
many infinitely unbounded solutions for $\gamma<0$. All these solutions are increasing and 
strictly concave and uniqueness of the bounded solution hold for $\gamma\leq 0$.

\item For $m\geq -\frac{1}{2}$ all the solutions are bounded.

\item For $-\frac{1}{2}\leq m \leq 1$ and for every $\gamma \in \mathbb{R}$ the problem 
(\ref{equation})-(\ref{cond03}) has one and only one solution. This solution is strictly concave
and increasing.

\item For $m>1$ and $\gamma \in \mathbb{R}$ the problem (\ref{equation})-(\ref{cond03}) 
has one and only one concave solution and infinitely many concave-convex solutions. 
Moreover, there is an unique concave-convex solution 
that verifies $f(t)\rightarrow \lambda>0$ as $t\rightarrow \infty$ and all the other 
concave-convex solutions are such that $f(t)\rightarrow 0$ as $t\rightarrow \infty$.
\bigskip
\end{itemize}
After this study it remains to investigate the following situations
\begin{itemize}
\item For $-1<m< -\frac{1}{2}$ and $\gamma >0$ is the bounded solution unique ?
\item For $-1<m< -\frac{1}{2}$ and $\gamma \geq 0$ is there unbounded solution ?
\end{itemize}

\end{document}